%% file: classification-monoperiodic.tex
\documentclass[american]{article}
\usepackage[utf8]{inputenc}
\usepackage{amsmath}
\usepackage{babel}
\usepackage{amsthm}
\usepackage{amssymb}
\usepackage{mathtools}
\usepackage{enumitem}
\usepackage{graphicx}
\usepackage{cutwin}
\usepackage{wrapfig}
\usepackage{hyperref}
\usepackage{cleveref}
\usepackage{caption}
\usepackage{changepage}
\usepackage[skip=1ex, labelfont=bf, font=footnotesize]{caption}
\usepackage{tikz}
\usepackage{makecell, tabularx}
\usepackage[backend=bibtex]{biblatex}
\usepackage{xcolor}
\usetikzlibrary{hobby, calc, math, intersections, decorations.markings, decorations.pathreplacing, arrows, angles, external}
\usepackage{tkz-euclide}
\input{insbox}

\setcellgapes{4pt}
\newcolumntype{Y}{>{\centering\arraybackslash}X}

\newcommand{\Z}{\mathbb{Z}}
\newcommand{\R}{\mathbb{R}}

\newcommand{\N}{\mathbb{N}}
\newcommand{\cupdiag}{\substack{\cup \\ \cap}}

\DeclareRobustCommand{\revPer}{\text{\reflectbox{$\%$}}}

\DeclareMathOperator{\spread}{sp}
\DeclareMathOperator{\len}{len}
\DeclareMathOperator{\id}{id}
\DeclareMathOperator{\MCG}{MCG}
\DeclareMathOperator{\au}{au}
\DeclareMathOperator{\lk}{lk}
\DeclareMathOperator{\sw}{sw}
\DeclareMathOperator{\red}{red}
\DeclareMathOperator{\conv}{conv}
\DeclareMathOperator{\sgn}{sgn}
\DeclareMathOperator{\wrap}{wr}

\newtheorem{theorem}{Theorem}[section]
\newtheorem{prop}[theorem]{Proposition}
\newtheorem{lemma}[theorem]{Lemma}
\newtheorem{corollary}[theorem]{Corollary}
\newtheorem{definition}[theorem]{Definition}
\newtheorem{conjecture}[theorem]{Conjecture}

\title{Classification of small links\\ in the unmarked solid torus}
\author{John M. Sullivan, Max Zahoransky von Worlik\\
\small Institut f\"ur Mathematik, Technische Universit\"at Berlin}
\date{12 June 2023}

\begin{document}

\maketitle
\begin{abstract}
	We introduce a framework to analyze knots and links in an unmarked solid torus. We discuss invariants that detect when such links are equivalent under an ambient homeomorphism, and show that the multivariable Alexander polynomial is such an invariant. We compute, for links with low wrapping number, bounds on the degree of a Dehn twist needed to transform one into the other that depend on the dichromatic Kauffman polynomial. Finally, we use this to give a classification of all non-split links up to 6 crossings in the unmarked solid torus.
\end{abstract}
 \section{Introduction}
 	A \emph{link} in a 3-manifold $M$ is a (smooth) embedding of a finite number of disjoint 1-spheres into $M$. A \emph{classical link} is a link in the manifold $S^3$. We are interested in links in the manifold $D^2 \times S^1$, the solid torus.
 	
 	Links in the solid torus are of particular interest as they are related to certain periodic structures in 3-dimensional space: Any link $L \subset D^2 \times S^1$ lifts to a link $\widetilde L$ in the universal cover $D^2 \times \R \subset \R^3$. In this way, one may consider the study of solid torus links as an extension of classical knot theory in $\R^3$ to certain unbounded knotted structures.
 	
 	There have been previous enumerations of knots \cite{gab12} and links \cite{sch16} in the solid torus. The key difference in our approach is that rather than classifying links in the marked solid torus, we classify them in the unmarked solid torus. In other words, we consider links to be equivalent if there is an ambient homeomorphism taking one to the other, whereas \cite{gab12, sch16} required an ambient isotopy. This means that we will identify links that differ only by a Dehn twist of the solid torus. We are also interested in listing composite links, whereas \cite{gab12, sch16} only gave lists of prime knots and links.
 	
 	In order to determine which links are equivalent in this sense, we need to use link invariants that can detect this equivalence. One well-established such invariant is the hyperbolic volume of the link complement \cite{thurnotes}, but it proves insufficient in many cases. In section 3, we show that the multivariable Alexander polynomial offers a way to distinguish links that are not homeomorphically equivalent. In section 4, using the two-variable Kauffman polynomial \cite{HP89}, we compute explicit bounds on the degree of a possible Dehn twist relating two links.
 	
 	To get our classification, we build on the algorithm presented in \cite{gab12}. In contrast to it, we use slightly different, periodic diagrams (and their associated periodic Gauss paragraphs) to generate all possible small links. These consist of curves that start and end in 2 boundary lines, as well as circular segments. We give some technical results on this representation in section 5.
 	
 	Then we start by classifying up to ambient isotopy, using mainly the dichromatic Kauffman polynomial, followed by determining which of the remaining links are homeomorphically equivalent using the methods established in sections 3 and 4. The algorithm we use to produce our link table is explained in detail in section 7. Finally, section 8 contains the resulting link table, with some of the corresponding invariants given in the appendix. \\
 	The code for the algorithm is written in Python and can be found online \cite{code}.
 	\subsection*{Acknowledgements}
 	This research was partially supported by the DFG Collaborative Research Center TRR 109 ``Discretization in Geometry and Dynamics''.
 	
 	\section{Links in the solid torus}
 	A \emph{$k$-component link} $L$ in an oriented 3-manifold $M$ is a closed smooth 1-dimensional submanifold with $k \in \N$ path components. If $M$ is a manifold with boundary, we require that $L$ is contained in the interior of $M$. If $L$ is equipped with an orientation, we call $L$ an \emph{oriented link}. A 1-component link is also called a \emph{knot}.
 	
 	We say that a link $L$ in $M$ is \emph{affine} if it is contained in a 3-ball. As in the case of classical links, we may form a \emph{connected sum} $L \# A$ of a link $L$ in $M$  and an affine link $A$, by cutting both links open and sewing $A$ into $L$. Note that we do not allow this operation for two non-affine links\footnote{One might wonder why we restrict ourselves to decompositions where one component is affine; the issue here is twofold. First, in order to build the connected sum of two links contained in some manifold $M$, we also have to attach the manifold to itself, and in general it is not clear that this can be done in a way such that we do not end up with a different manifold. Secondly, even where this is possible it leads to undesirable outcomes, e.g. in the solid torus it would be possible to obtain the affine unknot as the connected sum of two non-affine knots.}, and the operation depends on a choice of component of $L$ along which the surgery is performed. We say that $L$ is \emph{composite} if it can be written as such a connected sum, where $L,A$ are not affine unlinks, and \emph{prime} otherwise.
 	
 	In the following, we will write the solid torus as $V^3 := D^2 \times S^1$. For links in $S^3$, there are two different notions of (oriented) link equivalence that give rise to the same equivalence relation: One says that two links $K, L$ are equivalent if there is an orientation-preserving homeomorphism $h: S^3 \to S^3$ such that $h(K) = L$, or if there is an ambient isotopy of $S^3$ taking $K$ to $L$. In the solid torus $V^3 := D^2 \times S^1$ (with fixed orientation), these notions will yield different equivalence relations.
 	\begin{definition}
 		Let $K,L$ be two links in the solid torus.
 		\begin{enumerate}
 			\item $K$ and $L$ are \emph{homeomorphically equivalent} or \emph{$h^+$-equivalent} if there is an orientation-preserving homeomorphism $h: V^3 \to V^3$ such that $h(K) = L$.
 			\item $K$ and $L$ are \emph{ambiently equivalent} or \emph{a-equivalent} if there is an ambient isotopy $F: V^3 \times [0,1] \to V^3$ such that $F_0 = \id_{V^3}$ and $F_1(K) = L$.
 		\end{enumerate}
 	\end{definition}
 	Occasionally, in the literature \cite{dnpr18} the former is just called equivalent and the latter isotopic; we opted for different naming to avoid a double meaning of the term equivalent.
 	
 	Note that a-equivalence still implies h-equivalence. As is common in the classical case, we refer to both a concrete submanifold and its a-equivalence class as a link, only indicating the distinction when necessary.
 	
 	We also want to identify links which are mapped onto one another by orientation-reversing homeomorphisms. In practice this is also the equivalence class most often (if somewhat tacitly) used for classical knots: The standard knot tables list only the trefoil knot and drop its mirror image, even though the two are not a-equivalent. The utility of restricting oneself to orientation-prevesering homeomorpisms in the classical case lies in the fact that these are the homeomorphisms that can be realized via isotopies. In the solid torus we have already seen that these notions do not coincide, so we might as well allow our homeomorphisms to be orientation-reversing. We call two links related by any ambient homeomorphism \emph{h-equivalent}.
 	
 	To better enable us to talk about the homeomorphisms of $V^3$, we introduce the \emph{mapping class group} $\MCG(V^3)$, the group of orientation-preserving homeomorphisms of $V^3$ up to isotopy. A detailed introduction to mapping class groups in general can be found in \cite{fm12}. We will need the following theorem, proved in \cite{br98}:
 	\begin{theorem}
 		\label{mappingclassgroup}
 		The mapping class group of the solid torus is
 		\begin{equation*}
 			\MCG(V^3) \cong \Z_2 \times \Z.
 		\end{equation*}
 	\end{theorem}
    The generator of the torsion subgroup is the combination of two reflections in $D^2$ and $S^1$, respectively, and can be thought of as flipping a donut over onto its glazed side. We will call this homeomorphism $r$ throughout the paper.
    
    The generator of the free subgroup is the \emph{1-fold Dehn twist}, the result of cutting the solid torus open along a disk, twisting it around a full rotation, and then gluing the disks back together We will call this generator $d$.
    
	Mapping class groups are mostly studied for orientation-preserving maps, but of course we want to include all homeomorphisms. We call the group of all homeomorphisms the \emph{extended mapping class group} $\MCG^{\pm}(V^3)$. The standard mapping class group from before is an index 2 subgroup of this, as the composition of any two orientation-reversing maps will be orientation-preserving, so to get a generating set for $\MCG^{\pm}(V^3)$ it suffices to add one element of $\MCG^{\pm}(V^3) \setminus \MCG(V^3)$ to $\{r,d\}$; we choose the reflection in the $S^1$ factor of the solid torus, which we call $q$. It is not hard to see that we have the following presentation:
	\begin{equation}
		\MCG^{\pm}(V^3) \cong \left\langle q,r,d \, \middle| \, q^2, r^2, (qr)^2, (qd)^2, rdrd^{-1} \right\rangle \cong \Z_2 \times D_\infty
	\end{equation}
	It is worth noting that the impact of the torsion elements of $\MCG^\pm(V^3)$ is functionally quite different from the Dehn twist $d$; thus it will sometimes be convenient to distinguish them.
	\begin{definition}
		Two links in $V^3$ are \emph{symmetric} if they are related by an ambient homeomorphism contained in the subgroup of $\MCG^\pm(V^3)$ generated by $r$ and $q$.
	\end{definition}
	
 	A \emph{marking} of the solid torus is a longitude in the boundary torus, i.e. an oriented simple closed curve in $\partial V^3$ that represents a generator of $H_1(V^3)$. Two markings are equivalent if they are in the same homology class in $H_1(\partial V^3)$. A \emph{marked} solid torus is the solid torus together with a marking.
 	
 	For any two marked solid tori, there is an orientation-preserving homeomorphism that takes one marking to the other. If we disregard the orientation of the marking, this can always be achieved via a homeomorphism isotopic to $d^n$ for some $n$.
 	
 	As in the classical case, we want to study solid-torus links by means of link diagrams. Given a marked solid torus, we may embed it into $\R^3$ (with standard orientation) in such a way that the marking is mapped to $S^1 \times \{0\}$ oriented counterclockwise, the embedding respects the orientations of $V^3$ and $\R^3$, and the image of the embedding is a tube of radius $\tfrac 1 4$ around $\tfrac 3 4 S^1 \times \{ 0\}$. We call this embedded torus the \emph{round solid torus}. Any two such embeddings are isotopic.
 	
 	The image of a link $L$ in the round solid torus can be projected to an annulus $A$ in the $xy$-plane; the result is a graph $\Gamma$ embedded in $A$. By a small isotopy, we can ensure that the map $L \to A$ is a smooth immersion that is one-one except at the vertices, all vertices of $\Gamma$ are 4-valent, and at any vertex the angles between outgoing edges are nonzero. We call the vertices of $\Gamma$ \emph{crossings}. $\Gamma$ may be decorated at the crossings with over- and undercrossing information to give an \emph{annulus diagram} $D$ for $L$. Instead of drawing the annulus, we may simply mark the face of the diagram that contains the inner cavity of the annulus. We say that the planar graph $\Gamma$ is the \emph{underlying graph} of $D$. We note that $\Gamma$ may contain circles that do not coincide with any vertex, and thus calling it a graph may seem questionable; for our purposes it is fine to allow edges without any boundary in our graph.
 	
 	As in the classical case, Hoste and Przytycki have shown\cite{HP89} that two annulus diagrams belong to a-equivalent links if and only if they can be transformed into one another by a finite sequence of local isotopy moves, the \emph{Reidemeister moves}.
 	
 	In particular for dealing with Dehn twists, it will be convenient to consider another diagrammatic representation: In the annulus, pick any simple curve connecting the boundary components that does not intersect any crossings. We may cut open the annulus diagram along this curve, which yields (after homeomorphism) a diagram in a square, which we call a \emph{periodic diagram}. This name is chosen since the diagram can be thought of as part of a projection of the periodic structure we obtain by lifting the link to the universal cover $D^2 \times \R$.\\
 	If we start with an oriented link $L$, then the annulus and periodic diagrams inherit an orientation from $L$; we call a diagram with this additional information \emph{oriented}.
 	\begin{definition}
 		Let $D$ be an annulus diagram of a link $L$ with underlying graph $\Gamma$, and $P$ the periodic diagram obtained from $D$ by cutting along a curve $\gamma$ that is transverse to $D$ (thus does not meet any vertex).
 		\begin{itemize}
 			\item A \emph{passing} of $(D, \gamma)$ is an intersection point of $\Gamma$ and $\gamma$. A passing of $P$ is the corresponding pair of points in $P$ on the boundary of the square.
 			\item A \emph{segment} of $(D, \gamma)$ (or $P$) is either the image of an arc between two passings that contains no other passing, or the image of a link component whose image does not meet $\gamma$. In the latter case, we say that the segment is \emph{closed}.
 			\item The \emph{wrapping number} of $(D, \gamma)$ is its number of passings. The \emph{wrapping number} of $D$ is the minimum of the wrapping numbers of $(D, \widetilde{\gamma})$ for any splitting curve $\widetilde{\gamma}$. The \emph{wrapping number} of $L$ is the minimum wrapping number across all its diagrams.
 		\end{itemize}
 	\end{definition}
 	The two simplest knots in $V^3$ we can imagine are the \emph{(affine) unknot}, the unique affine knot that has a diagram without crossings, and the \emph{longitude knot}, given by the embedding $f: S^1 \to D^2 \times S^1 = V^3, f(x) = (0,x)$. They have standard periodic diagrams with no crossings and exactly one segment: For the unknot, this is the diagram $D_\circ$ consisting of just a circle not touching the boundary, and for the longitude knot the diagram $D_|$ consisting of a single vertical line.
 	
 	We will point out here that links with small wrapping number fall into simple categories: Links with wrapping number 0 are exactly the affine links, and links with wrapping number 1 always have a 2-sphere intersecting them twice that contains all the knotting, or more precisely, they are connect sums of the longitude knot with some affine link.
 	
 	Note also that we may specify an orientation for a periodic diagram by giving each segment some orientation. An orientation of the annulus diagram $D$ will always induce an orientation on the periodic diagram $P$, but vice versa not all possible orientations of $P$ will induce a valid orientation of the link.
 	\begin{definition}
 		An annulus diagram of a link is \emph{split} if its underlying graph is disconnected. A link is \emph{split} if it has a split annulus diagram. A link is \emph{h-split} if it is h-equivalent to a split link.
 	\end{definition}
 	The components of a split link can be separated to lie in two distinct smaller solid tori in $V^3$ which are not themselves knotted or linked.
 	\begin{definition}
 		The \emph{crossing number} of an annulus diagram is the number of vertices in its underlying graph. The \emph{crossing number} of a link is the minimum crossing number of any annulus diagram of the link.
 	\end{definition}
 	\begin{theorem}
 		Two periodic diagrams $F$ and $F'$ give rise to a-equivalent links if and only if they can be transformed into one another by a series of the following moves:
 		\begin{center}
 			\input{tikz/omega-theta-move.tex}
 		\end{center}
 	\end{theorem}
 \begin{proof}
 	We consider the annulus diagrams $\tilde{F}, \tilde{F}'$ obtained by gluing $F$ rsp. $F'$ together along the top and bottom. If the corresponding links are a-equivalent, then there is a sequence for Reidemeister moves $\Omega_i$ and isotopies of $A$ taking $\tilde{F}$ to $\tilde{F}'$. If the isotopies change the position of the diagram relative to the gluing line, we obtain moves of the form $\Theta_i^\pm$. Conversely, if there is a sequence of moves $\Omega_i, \Theta_i^\pm$, then these give a Reidemeister sequence from $\tilde{F}$ to $\tilde{F}'$.
 \end{proof}
\begin{corollary}
	\label{hequivmoves}
	Two periodic diagrams $F$ and $F'$ give rise to h-equivalent links if and only if they can be transformed into one another by a series of moves of type $\Omega_i, \Theta_i^+$ as above together with
the $180^{\circ}$ rotation $R$ of the diagram,
the reflection $Q$ of the diagram across a horizontal line,
and the move $\Delta$ shown below.
	\begin{center}
		\input{tikz/delta-move.tex}
	\end{center}
\end{corollary}
\begin{proof}
	First assume that we have a sequence of such moves. The moves $\Omega_i, \Theta_i$ are isotopy moves. $\Delta$ corresponds to the Dehn twist homeomorphism $d$, whereas the $180^\circ$ rotation $R$ corresponds to $r$, and finally the reflection $Q$ corresponds to a reflection in the $S^1$ factor, i.e. $q$. All of these are ambient homeomorphisms, so every step of the process is an h-equivalence move.\\
	The other direction follows since $\{ d, r, q \}$ generates $\MCG^\pm(D^2 \times S^1)$, and we have a diagram move corresponding to each of these. We may drop the moves $\Theta_i^-$ since we have that $\Theta_1^- = Q \circ \Theta_1^+ \circ Q$ and $\Theta_2^- = R \circ \Theta_2^+ \circ R$.
\end{proof}

We can see that the wrapping number of a link $L$ is in fact an invariant of h-equivalence: This is since rotation, reflection and the move $\Delta$ from above leave the wrapping number of a diagram invariant.

We will finally note that any link has a special type of diagram, a fact that will become critical when trying to generate all diagrams up to a given crossing number.

\begin{definition}
		A \emph{basic} diagram $D$ is a periodic diagram such that any two segments of $D$ cross at most once, and there are no self-crossings.
\end{definition}

Note that a closed segment intersects each other segment an even number of times. Thus a closed segment in a basic diagram represents an unknotted affine component split from the rest of the link.

\begin{prop}
	\label{passBound}
	Let $L$ be a non-split link in the solid torus with crossing number $c$ and wrapping number $\wrap$. Then there is a basic diagram of $L$ whose number of passings is at most $2c + 2 - \wrap$, and thus $\wrap \leq c + 1$.
\end{prop}
\begin{proof}
	Consider an annulus diagram $D$ of $L$ realizing the crossing number with underlying embedded planar graph $\Gamma$. We know that we get a periodic diagram from this by specifying a simple curve from the marked face to the unbounded face. We can find such a curve $\gamma$ in the plane crossing a minimal number of edges, and $\gamma$ will correspond to a path in the dual planar graph $\Gamma^*$, which we will also call $\widetilde \gamma$.
	
	Note that $\gamma$ crosses each face at most once; otherwise we could reduce the number of edges crossed by $\gamma$, replacing a path from a face to itself with a path inside the face. Thus the induced path $\widetilde \gamma$ is simple. We write $\len \gamma$ for the number of edges crossed by $\gamma$; it is immediate that $\len \gamma \geq \wrap$. Furthermore, as $L$ is non-split, $\Gamma$ is a connected planar graph with $c$ 4-valent vertices, and hence by Euler's formula has $c+2$ faces. Since $\gamma$ meets each face at most once, we have $\len \gamma \leq c+1$.
	
	Since $\widetilde \gamma$ is a simple path, we may choose a spanning tree $T$ for $\Gamma^*$ that contains $\widetilde \gamma$; this graph may be embedded in the annulus diagram $D$, with a vertex in each corresponding face, in such a way that $\widetilde \gamma \subseteq T$ is mapped back to $\gamma$.
	\begin{figure}
		\begin{center}
			\input{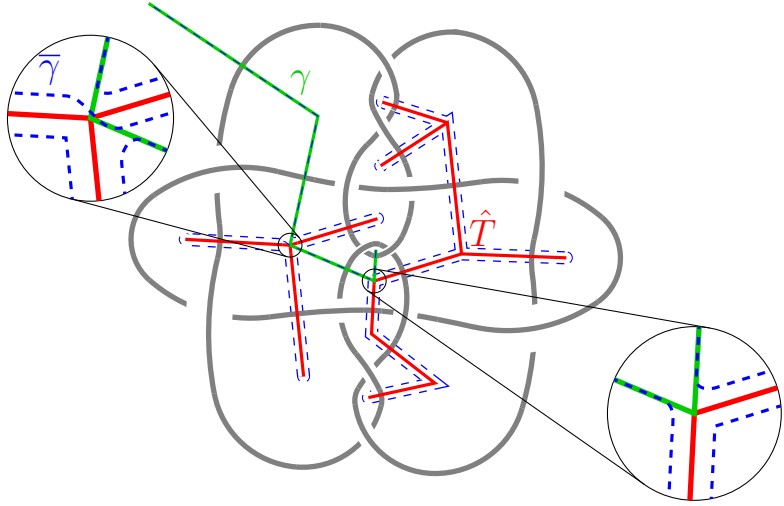}
			\caption{The new cutting curve $\overline{\gamma}$ constructed from $T$. The marked region of the annulus diagram is the central bigon.}
			\label{branchGamma}
		\end{center}
	\end{figure}
	We will now construct a new cutting curve. For this, we say that a vertex $v$ of $T$ is \emph{branching} if $v \in \widetilde \gamma$ and $T \setminus \widetilde \gamma \neq \emptyset$. Let $N$ be a sufficiently small neighborhood of $T \setminus \gamma$ in $\R^2$; then $\partial N$ forms a simple curve away from branching vertices.
	
	The new cutting curve $\overline \gamma$ is constructed by joining $\partial N$ with $\gamma$ along branching vertices: We follow $\gamma$ until we reach a branching vertex, then we make a small arc to the part of $\partial N$ to the left of $\widetilde \gamma$ and follow $\partial N$ along all branches of $T$ that lie to the left of $\gamma$. Upon returning, we switch to the right and follow $\partial N$ along all brnaches to the right; once all branches emerging from the branchng vertex have been visited, we continue along $\gamma$. This can be done in a way that leaves the resulting curve simple, see $\cref{branchGamma}$.
	
	Now it is not hard to see that the curve $\overline{\gamma}$ defined by this process is a cutting curve that will yield a basic diagram: For any non-basic diagram has two segments that intersect twice, and hence enclose a face, which corresponds to a face in the annulus diagram. But $\overline{\gamma}$ intersects all faces, so we cannot have an untouched face in the interior of the diagram.
	
	As $\Gamma^*$ has $c+2$ vertices, so does $T$ and as a tree it then has $c+1$ edges. $\overline \gamma$ runs twice along each edge of $T$ except those contained in $\gamma$, which it only passes once. Thus 
	\begin{equation} 
		\len \overline{\gamma} = 2 (c+1) - \len \gamma  \leq 2c + 2 - \wrap.
	\end{equation}
	But we clearly have that $\wrap \leq \overline \gamma$ always, thus from the above we may conclude that $\wrap \leq 2c + 2 - \wrap$, and hence $\wrap \leq c+1$.
\end{proof}

 \section{The Alexander polynomial for solid torus links}
 We can augment any link $L$ in a marked solid torus to obtain a colored link in $S^3$, by embedding $L$ via the marking embedding $g: V^3 \to \R^3 \subset S^3$ into the round solid torus and adding another component $m$ corresponding to the z-axis extended to $\infty$. The \emph{augmented link} $\au(L)$ is the link $g(L) \cup m \subseteq S^3$, with $m$ and $g(L)$ having different colors. If $L$ is an oriented link, we need an orientation for $m$ to get an oriented link in $S^3$; we choose the orientation of the z-axis going up. Note that $S^3$ inherits an orientation from the standard orientation of $\R^3$. Clearly $\au(L)$ is an invariant of a-equivalence.
 
 The maps $q, r: V^3 \to V^3$ from \cref{mappingclassgroup} then extend to maps $\overline q, \overline r: S^3 \to S^3$ in the sense that $g \circ q = \overline q \circ g$ and $g \circ r = \overline r \circ g$; here $\overline r$ is a rotation in $S^3$ and $\overline q$ a reflection. Then $\overline r$ is isotopic to the identity and thus for an unoriented link $L \subset V^3$, we have that $\au(r(L)) = \au(L)$. If we add orientation, note that the only thing that changes is the orientation of the component $m$; i.e. $\au(r(L))$ is the same as $\au(L)$ with the orientation on $m$ reversed.
 
 Analogously, $\au(q(L)) = \overline q(\au(L))$; here the orientation on $m$ does not change. $\overline q$ is orientation-reversing and thus generically $\overline q(\au(L))$ will not be isotopic to $\au(L)$.
 
 It is also immediate that the homeomorphism type of the complement $V^3 \setminus L$ is an invariant of h-equivalence. At the same time, $V^3 \setminus L \cong S^3 \setminus \au(L)$. Thus, if the complement admits a hyperbolic structure, we can use the theory of hyperbolic links \cite{thurnotes} to derive many invariants of h-equivalence, such as the hyperbolic volume. However, we should be mindful that in this process, the coloring information is lost.
 
 Next we will use the \emph{linking number} of two oriented knots in $S^3$: Let $K_1, K_2 \subset S^3$ be disjoint knots. It is well known that $H_1(S^3 \setminus K_1) \cong \Z$, and the group is generated by a meridian winding counterclockwise around the link (i.e. given a frame for $K_1$ compatible with the orientation of $S^3$, the meridian winds counterclockwise in the normal plane with induced orientation). Then $K_2$ generates a homology class in $H_1(S^3 \setminus K_1)$, so $[K_2] = k\beta$ for some $k \in \Z$. We call $k$ the \emph{linking number} $\lk(K_1, K_2)$. This generalizes to the case where $K_2$ is a multi-component link.
 
 Using the linking number allows us obtain new invariants of h-equivalence.
 \begin{definition}
 	Let $L \subset V^3$ be an oriented link with $k$ components $L_1, \ldots, L_k$, and $\au(L)$ its augmentation by  the additional component $m$. The \emph{winding number} $w_j$ of $L_j$ is the linking number $\lk(m, L_j)$. The \emph{total winding number} of $L$ is $w(L) = \lk(m,L)$.
 \end{definition}
\begin{prop}
	The unordered collection $w_1, \ldots, w_k$ of winding numbers, as well as $w(L)$, are invariants of $h^+$-equivalence.
\end{prop}
 \begin{proof}
 	The winding number of a knot $L_i$ in $V^3$ is its homology class $[L_i] \in H_1(V^3) \cong \Z$. This is clearly invariant under isotopy, and to see that it is also invariant under h-equivalence it suffices to point out that any homeomorphism $h: V^3 \to  V^3$ induces an isomorphism on $H_1(V^3)$, which is the identity if $h$ is orientation-preserving.
 \end{proof}
 \begin{prop}
 	\label{windingWrapping}
 	Let $L \subset V^3$ be an oriented link with wrapping number $\wrap$. Then $w(L) \leq \wrap$.
 \end{prop}
 \begin{proof}
 	Consider an annular diagram $D$ of $L$ and a cutting curve $\gamma$ realizing the wrapping number as the geometric intersection number of $D$ with $\gamma$.  It is well known that the linking number $\lk(m,L)$ equals the algebraic intersection number in $S^3$ of $L$ with a meridian disk bounded by $m$.  Projecting to the annulus, that is the algebraic intersection number of the oriented diagram $D$ with any cutting curve, in particular with $\gamma$. This algebraic intersection number is at most the geometric intersection number.
 \end{proof}
 This invariant captures the coloring, but it is still very coarse, and distinct links will frequently have the same winding numbers. For a finer invariant, we turn to the multivariable Alexander polynomial $\tilde{\Delta}_L$ of a classical oriented link, a polynomial with one variable for each link component.
 
 The Alexander polynomial is defined as follows: For an oriented $k$-component link $L$, we may choose generators $t_1, \ldots, t_k$ of $H_1(S^3 \setminus L) \cong \Z^k$, each a meridian around one of the components with linking number $+1$. Note that the Hurewicz map $h: \pi_1(S^3 \setminus L) \to H_1(S^3 \setminus L)$ is onto, and $\ker h$ corresponds to a covering $p: X \to S^3 \setminus L$ with $p_*(\pi_1(X)) = \ker h$. Then $H_1(S^3 \setminus L)$ acts on $X$ via deck transformations, and thus provides a scalar multiplication on $H_1(X)$, turning $H_1(X)$ into a module over the integral group ring $\Z H_1(S^3 \setminus L)$. This ring can be understood as the ring $\Lambda$ of integer Laurent polynomials in the variables $t_i$. Also the relative homology group $H_1(X, p^{-1}(b_0))$, where $b_0 \in S^3 \setminus L$ is some base point, becomes a $\Lambda$-module in the same way.
 
 Any $\Lambda$-module $H$ has a free resolution, i.e. an exact sequence
 \begin{equation}
 	\label{eqExact}
 	\Lambda^n \stackrel{P}{\to} \Lambda^m \to H \to 0
 \end{equation}
with free $\Lambda$-modules. The map $P$ may be viewed as an $n \times m$ matrix with $\Lambda$ coefficients. Choosing such a free resolution for $H_1(X, p^{-1}(b_0))$, the multivariable Alexander polynomial $\tilde \Delta_L(t_1, \ldots, t_k) \in \Lambda$ of $L$ is defined as the greatest common divisor of the determinants of all maximal square submatrices of $P$. This turns out to be an invariant of $L$ up to multiplication by factors of the form $\pm \prod_i t_i^{s_i}$; for a proof and detailed exposition we refer the reader to \cite[Chapter 7]{ka12}. We write $p_1 \doteq p_2$ if 2 polynomials are equal up to such a factor.
\begin{prop}
	\label{oriRev}
	Let $L \subset S^3$ be an oriented link with component knots $L_1, \ldots, L_k$, and let $L^{-i}$ denote that link with the orientation of $L_i$ reversed. Then
	\begin{equation*}
		\tilde \Delta_{L^{-i}}(t_1, \ldots, t_k) \doteq \tilde \Delta_L(t_1, \ldots, t_{i-1},t_i^{-1},t_{i+1}, \ldots, t_k).
	\end{equation*}
\end{prop}
\begin{proof}
	This follows from our definition of the Alexander polynomial: Reversing the orientation on $L_i$ means that, instead of the basis $t_1, \ldots, t_k$ of meridians for $H_1(S^3 \setminus L)$, we choose a differently oriented meridian for $L_i$, thus our basis becomes $t_1, \ldots, t_{i-1},t_i^{-1},t_{i+1}, \ldots, t_k$.
\end{proof}
We will make use of a classical result first proven by Torres \cite{to53}.
\begin{theorem}[First Torres condition]
	\label{torres1}
	Let $L \subset S^3$ be an oriented link, then
	\begin{equation}
		\tilde{\Delta}_L(t_1, \ldots, t_k) \doteq \tilde{\Delta}_L(t_1^{-1}, \ldots, t_k^{-1}).
	\end{equation}
\end{theorem}
\begin{theorem}
	\label{alexSym}
	Let $L \subset S^3$ be an oriented link, $m$ a reflection in $S^3$, and $-L$ be the link with reversed orientation. We have that
	\begin{equation}
		\tilde \Delta_L(t_1, \ldots, t_k) \doteq \tilde \Delta_{-L}(t_1, \ldots, t_k) \doteq \tilde \Delta_{m(L)}(t_1, \ldots, t_k).
	\end{equation}
\end{theorem}
\begin{proof}
	Reversing the orientation of $L$ is the same as reversing each component of $L$; thus by \cref{oriRev} it follows that
	\begin{equation}
		\tilde \Delta_{-L}(t_1, \ldots, t_k) \doteq \tilde \Delta_L (t_1^{-1}, \ldots, t_k^{-1}),
	\end{equation}
	and together with \cref{torres1}, the first equation follows.
	
	For the second equation, note that $m$ induces the map
	\begin{equation}
		m_*: \pi_1(S^3 \setminus L) \to \pi_1(S^3 \setminus m(L)), \qquad t_i \mapsto t_i^{-1}.
	\end{equation}
	of the fundamental groups, as $m$ sends meridians of $L$ to meridians of $m(L)$, but changes their linking number. As before, this implies that
	\begin{equation}
		\tilde \Delta_{L}(t_1, \ldots, t_k) \doteq \tilde \Delta_{m(L)} (t_1^{-1}, \ldots, t_k^{-1}).
	\end{equation}
	The claim again follows due to \cref{torres1}.
\end{proof}
We will use a two-variable version of the Alexander polynomial. Given a solid torus link $L$, the \emph{dichromatic Alexander polynomial} is
 \begin{equation}
 	\Delta_L(t,x) := \tilde{\Delta}_{\au(L)} (t,\ldots, t, x),
 \end{equation}
where we understand that the final variable of $\tilde{\Delta}_{\au(L)}$ corresponds to the additional component $m$. As the multivariable Alexander polynomial is only defined up to multiplication by a monomial with coefficient $\pm 1$, so too is $\Delta_L$.
\begin{prop}
	\label{revAlex}
	Let $L$ be an oriented solid torus link. Then we have for the dichromatic Alexander polynomial:
	\begin{equation}
		\Delta_L(t,x) \doteq \Delta_{q(L)}(t,x) \doteq \Delta_{r(L)}(t^{-1}, x)
	\end{equation}
\end{prop}
\begin{proof}
	As seen in the beginning of this section, $\au(q(L))$ is the mirror image of $\au(L)$, and by \cref{alexSym}, mirror images do not change the Alexander polynomial, so the dichromatic Alexander polynomial will remain the same as well.
	
	There we also saw that $\au(r(L))$ is the same as $\au(L)$ with the orientation on $m$ reversed. By \cref{alexSym} and \cref{oriRev}, this implies
	\begin{equation}
		\Delta_{r(L)}(t,x) \doteq \tilde \Delta_{\au(r(L))}(t, \ldots, t,x) \doteq \tilde \Delta_{\au(L)}(t, \ldots, t, x^{-1}) \doteq \Delta_L(t^{-1},x)
	\end{equation}
	and thus the claim.
\end{proof}
We now state the main theorem of this section, which gives the behaviour of $\Delta_L$ under the final generator of the mapping class group, $d$:
\begin{theorem}
 	\label{bicAlex}
 	Let $L$ be an oriented link in $V^3$. Then
 	\begin{equation}
 	\Delta_{d(L)} (t, x) \doteq \Delta_{L} (t, x t^w),
 	\end{equation}
 	where $w$ is the total winding number of $L$ and $d$ is the 1-fold Dehn twist.
 \end{theorem}
\begin{figure}
	\begin{center}
		\input{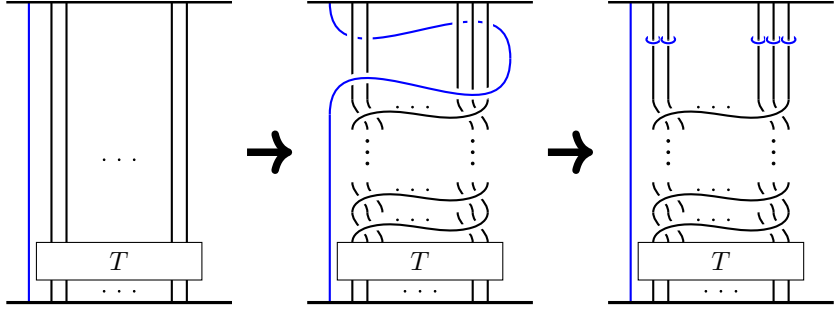}
	\end{center}
	\caption{The map $d$ transforms a meridian around $m$ into a curve homologous to $y\prod_{i=1}^k s_i^{w_i}$.}
	\label{figBicAlexProof}
\end{figure}
 \begin{proof}
 	We will proceed in three steps. First, we shall compute the induced isomorphism $d_*$ in homology. Then, we will argue that one may obtain a free resolution of the Alexander module of $\au(d(K))$ from the free resolution of that of $\au(K)$. Finally, we will combine these two results to prove our claim.
 	
 	\paragraph{Step 1} Let $k$ be the number of components of $L$, and $E := V^3 \setminus L \cong S^3 \setminus \au(L)$. Further let $E_d := V^3 \setminus d(L)$. There is a basis $t_1, \ldots, t_k, x$ of positive meridians for $H_1(E)$ and another one $s_1, \ldots, s_k, y$ for $H_1(E_d)$. It is evident that $d$ maps meridians of $L$ to meridians; i.e. we have $d_*(t_i) = s_i$ with appropriate indexing.
 	
 	Let $\gamma \subset V^3$ be a representative of $x$; the same curve is also a representative of $y$; as they are contained in $V^3$ it makes sense to draw them in a periodic diagram of $L$ rsp. $d(L)$ as in \cref{figBicAlexProof}. Now consider the image $d(\gamma)$. This will be $\gamma$ followed by a twist around $m$, that is, we can write (with multiplicative notation) $d_*(x) = y[\tilde{m}]$, where $\tilde{m}$ is a parallel of $m$. Then $[\tilde{m}] = y^{p_y}s_1^{p_1} \ldots s_k^{p_k}$ for some $p_y, p_1, \ldots p_k \in \Z$. Let $L_1, \ldots , L_r$ be the components of $L$. We have chosen our basis in such a way that the map $\varphi_i : H_1(E) \to H_1(S^3 \setminus L_j) \cong \Z$ induced by the embedding is just $\varphi_{i*}(t_j) = \delta_{ij}$ (the Kronecker delta) on the basis elements, and similarly $\varphi_{i*}(x) = 0$. 
 	
 	This implies that $\varphi_i ([\tilde{m}]) = p_i$, but on the other hand, this is just the homology class of $m$ in $S^3 \setminus L_i$, that is, $p_i$ is the linking number of $m$ and $L_i$, i.e. $w_i$. Similarly we can conclude that $p_y = 0$. Hence $d_*(x) = y \prod s_i^{w_i}$.
 	
 	\paragraph{Step 2} In order to compute the Alexander polynomial, we consider the universal abelian coverings $\widetilde E \stackrel{p}{\to} E$, $\widetilde{E_d} \stackrel{p_d}{\to} E_d$. We have the following picture:
 	\begin{center}
 		\begin{tikzpicture}
 			\node (eg) at (0,2) {$\widetilde E$};
 			\node (e) at (0,0) {$E$};
 			\node (te) at (2,0) {$E_d$};
 			\node (teg) at (2,2) {$\widetilde{E_d}$};
 			
 			\draw [->](eg) -- (e) node [midway, left] {$p$};
 			\draw [->] (teg) -- (te) node [midway, right] {$p_d$};
 			\draw [->] (e) -- (te) node [midway, below] {$d$}; 
 			
 			\coordinate (G1') at ($(eg.north west)$);
 			\coordinate (G2') at ($(teg.north east)$);
 			
 			\draw [->, shorten >=-1.5pt] (G1') arc (245:-70:1ex);
 			\draw [->, shorten >=-1.5pt] (G2') arc (-65:250:1ex);
 			
 			\node at ($(G1') + (-2em, 1em)$) (he) {$H_1(E)$};
 			\node at ($(G2') + (2em, 1em)$) (hte) {$H_1\left(E_d\right)$};
 			
 			\draw [->] (he) to [bend left=30] (hte);
 			\node at (1,3.7) {$d_*$};
 		\end{tikzpicture}
 	\end{center}
 	The deck transformation group action of $H_1(E)$ on $\widetilde E$ induces an action on $H_1\big(\widetilde E, p^{-1}(b_0)\big)$. This action extends to an action of the integral group ring $\Z H_1(E)$, and analogously this holds for $E_d$.
 	
 	Now $H_1(E)$ is a free abelian group of rank $k+1$. The integral group ring is then isomorphic to the integral Laurent polynomial ring $\Lambda = \Lambda\left[r_1, \ldots, r_k, z\right]$, with an isomorphism $g: \Lambda \to \Z H_1(E)$ induced by $g(z) = x, g(r_i) = t_i$. In the same way, $g_d: \Lambda \to \Z H_1\left( E_d \right)$ determined by $g_d(z) = y, \widetilde{g}(r_i) = s_i$ is an isomorphism.
 	
 	On the other hand, $d_*$ also extends to an isomorphism $\hat{d}_*: \Z H_1(E) \to \Z H_1\left( E_d \right)$. Let $\star: \Z H_1\left( E_d \right) \times H_1 \big( \widetilde{E_d}, p^{-1}(b_0) \big) \to H_1 \big( \widetilde{E_d}, p^{-1}(b_0) \big)$ denote the ring action, then we have two distinct ring actions of $\Lambda$ on $\widetilde{H} := H_1 \big( \widetilde{E_d}, p^{-1}(b_0) \big)$:
 	\begin{alignat}{2}
 		\star_1: \Lambda \times \widetilde{H} & \to \widetilde{H} & \star_2: \Lambda \times \widetilde{H} & \to \widetilde{H} \\
 		(p, x) & \mapsto g_d(p) \star x & (p,x) & \mapsto \left( \hat{d}_* \circ g \right) (p) \star x.
 	\end{alignat}
 	With each of these, $\widetilde{H}$ becomes a $\Lambda$-module. Let $\widetilde{H}_i$ denote the $\Lambda$-module obtained via scalar multiplication by $\star_i$. As in \cref{eqExact}, $\widetilde{H}_2$ has a free resolution
 	\begin{equation}
 		\Lambda^n \stackrel{P}{\to} \Lambda^m \stackrel{f}{\to}  \widetilde{H}_2 \to 0
 	\end{equation}
 	with presentation matrix $P = \left( p_{ij} \right)$, where the $p_{ij}$ are integral Laurent polynomials in the variables $z,r_1, \ldots, r_k$. Note that there is a ring automorphism on $\Lambda$ given by $a := g_d^{-1} \circ \hat{d}_* \circ g$. This in turn gives us a group isomorphism $a^{\times n}: \Lambda^n \to \Lambda^n$ via
 	\begin{equation}
 		a^{\times n}\left(p_1, \ldots, p_n\right) = \big( a \left(p_1 \right), \ldots, a \left(p_n \right) \big).
 	\end{equation}
 	Then, viewing $P$ as a group homomorphism (forgetting the module structure), there are unique maps $\widetilde{P}, \widetilde{f}$ such that the following diagram commutes:
 	\begin{center}
 		\label{commdiag}
 		\begin{tikzpicture}
 			\node (ln1) at (0,0) {$\Lambda^n$};
 			\node (ln2) at (0,2) {$\Lambda^n$};
 			\node (lm1) at (2,0) {$\Lambda^m$};
 			\node (lm2) at (2,2) {$\Lambda^m$};
 			\node (he1) at (4,0) {$\widetilde{H}$};
 			\node (he2) at (4,2) {$\widetilde{H}$};
 			\node (o1) at (5,0) {$0$};
 			\node (o2) at (5,2) {$0$};
 			\draw [->] (ln1) -- (lm1) node[midway, below] {$\widetilde{P}$};
 			\draw [->] (ln2) -- (lm2) node[midway, above] {$P$};
 			\draw [->] (ln2) -- (ln1) node[midway, left] {$a^{\times n}$};
 			\draw [->] (lm2) -- (lm1) node[midway, left] {$a^{\times m}$};
 			\draw [->] (lm1) -- (he1) node[midway, below] {$\widetilde{f}$};
 			\draw [->] (lm2) -- (he2) node[midway, above] {$f$};
 			\draw [->] (he2) -- (he1) node[midway, right] {$\id$};
 			\draw [->] (he1) -- (o1);
 			\draw [->] (he2) -- (o2);
 		\end{tikzpicture}
 	\end{center}
 	Since the top row is exact, so is the bottom row. Our claim is that $\widetilde{P}, \widetilde{f}$ are actually $\Lambda$-module homomorphisms, with $\widetilde{H}$ having module structure via $\star_1$.\\
 	It is a straightforward calculation to check that $\widetilde{P}$ actually has a matrix presentation via $\widetilde{P} = \left( a(p_{ij}) \right)$, so it is a module homomorphism. Let $q \in \Lambda$ and $v \in \Lambda^m$, then we have for $\widetilde{f}$
 	\begin{multline}
 		\widetilde{f} (q v ) = f \left( \left( a^{\times m} \right)^{-1} (q v) \right) = f \big( a^{-1}(q) \left( a^{\times m} \right)^{-1} (v) \big) \\ = \left( a^{-1}(q) \right) \star_2 \widetilde{f} (v) = q \star_1 \widetilde{f} (v),
 	\end{multline}
 	where the final equality holds since $a^{-1} = \left( \hat{d}_* \circ g \right)^{-1} \circ g_d$.
 	
 	\paragraph{Step 3} We have seen that
 	\begin{equation}
 		\Lambda^n \stackrel{\widetilde{P}}{\to} \Lambda^m \stackrel{\widetilde{f}}{\to} \widetilde{H}_1 \to 0
 	\end{equation}
 	is a free resolution of the $\Lambda$-module $\widetilde{H}_1$. For any square submatrix $P_\square$ of $P$ it follows that $\det\big( \widetilde{P}_\square \big) = a \det(P\square)$, where $\widetilde P_\square$ is the submatrix of $\widetilde P$ with the same pattern. But this means that the Alexander polynomial may be computed as
 	\begin{equation}
 		\begin{aligned}
 		\tilde \Delta_{\au(d(L))}(r_1, \ldots, r_k,z) & \doteq \gcd_{\square \mbox{ \footnotesize\,maximal}} \bigg(\det \big( \widetilde{P}_\square \big) \bigg) = \gcd_{\square \mbox{\footnotesize\,maximal}} \big( a \left( \det \left( P_\square \right) \right) \big) \\ & = a \bigg( \gcd_{\square \mbox{\footnotesize\,maximal}} \big( \det (P_\square) \big) \bigg) \doteq a \left( \tilde \Delta_{\au(L)}(r_1, \ldots, r_k,z)\right) \\ & = \tilde \Delta_{\au(L)} \left(r_1, \ldots, r_k, z\prod_{i=1}^k r_i^{w_i} \right).
 		\end{aligned}
 	\end{equation}
 	Here the last equation follows from the computation of $d_*$ in step 1. Setting $z=x, r_i=t$ yields the desired result.
\end{proof}

\section{Normalizing the Alexander polynomial}
Our goal is to extract an invariant of $h$-equivalence from the Alexander polynomial. We start by introducing some helpful concepts for dealing with Laurent polynomials. First, given a Laurent polynomial $p(x_1, x_2)$, recall that its \emph{$x_i$-degree} $\deg_{x_i}(p)$ is the largest power of $x_i$ that appears in the polynomial (which may be negative), or $\deg_{x_i}(p) = -\infty$ if $p=0$.

The \emph{total degree} of $p$ is then
\begin{equation}
	\deg(p) := \max \left\{ n_1 + n_2 \, \middle| \, x_1^{n_1} x_2^{n_2} \mbox{ is a monomial with nonzero coefficient} \right\},
\end{equation}
and as before $\deg(p) = -\infty$ if $p=0$.

In what follows, we will also need an ordering of integer pairs $[n] = (n_1, n_2)$. For this, we define the \emph{graded lexicographic order} as follows: We say $[n] < [m]$ if either $n_1+n_2 < m_1 + m_2$ or $n_1 + n_2 = m_1 + m_2$ and $n_1 < m_1$. With this we can give an order of integer polynomials:
\begin{definition}
	Let $p, q \in \Z[t,x]$.
	\begin{enumerate}
		\item The \emph{coefficient tuple} of
		\begin{equation}
			p(t,x) = \sum_{m = 0}^{\deg(p)}\sum_{n_t + n_x = m} a_{n_t, n_x}t^{n_t}x^{n_x}
		\end{equation}
		is the tuple of $a_{n_t, n_x}$ (including coefficients equal to zero) sorted by their index in graded lexicographic order.
		\item We say that $p < q$ if either $\deg p < \deg q$, or $\deg p = \deg q$ and the coefficient tuple of $p$ is lexicographically smaller than the coefficient tuple of $q$.
	\end{enumerate}
\end{definition}

Form this, we obtain a preliminary normalization for a Laurent polynomial.

\begin{definition}
	Let $p(t,x)$ be an integer Laurent polynomial. Let $S =  \{ \tilde p \in \Z[t,x] \, | \, \tilde p \doteq p \}$. Then the minimal element $p^N$ of $S$ is the \emph{normal form} of $p$.
\end{definition}
Notice the requirement that $p^N \in \Z[t,x]$, i.e. there are no negative powers in $p^N$. The set of possible degrees $\{ - \infty \} \cup \N$ is well-ordered, so there is an element in $S$ of minimal degree, and there can be only finitely many elements of any given degree in $S$. This guarantees the existence of $p^N$.

All these definitions apply in the case where $p$ is an Alexander polynomial. Now, let $L$ be an oriented solid torus link, let
\begin{equation}
	D\Delta(L):= \left\{ \Delta_{d^k(L)}^N(t,x) \, \middle| \, k \in \Z \right\},
\end{equation}
where $d$ is the Dehn twist.
\begin{prop}
	$D\Delta(L)$ contains a minimal element.
\end{prop}
\begin{figure}
	\begin{center}
	\input{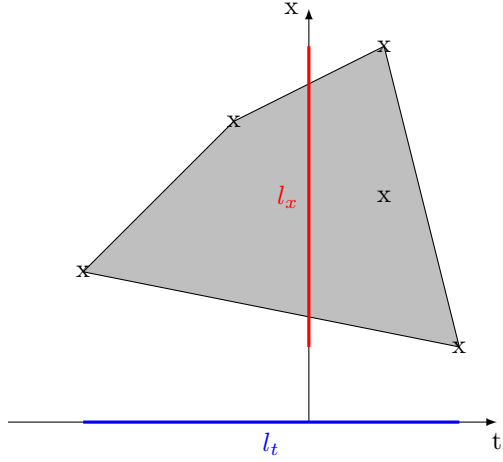}
	\caption{The Newton polytope of the polynomial $p(t,x) = 5t^2x - tx^5  + tx^3 - t^{-1}x^4 + 2t^{-3}x^2$.}
	\label{figNewtonP}
	\end{center}
\end{figure}
\begin{proof}
	Note first that we may assume a winding number $w(L) \geq 0$; otherwise we use the reflection $q$, which changes the sign of the winding number. Since $d \circ q = q \circ d^{-1}$, and by \cref{alexSym} we have that $\Delta_L(t,x) \doteq \Delta_{q(L)}(t,x)$, it follows that $\Delta_{d^k \circ q(L))}(t,x) \doteq \Delta_{d^{-k}(L)}(t,x)$, thus we get that $D\Delta(L) \doteq D\Delta(q(L))$. If $w=0$ or $\deg_x(\Delta_L^N(t,x)) \leq 0$ (i.e. there are no $x$ terms), then \cref{bicAlex} implies that $D\Delta(L)$ is actually a singleton set.
	
	Otherwise, consider for any polynomial $p \in \Lambda (t,x)$ the set
	\begin{equation*}
		\mathcal N(p) := \left\{ (n_t, n_x) \, \middle| \, t^{n_t}x^{n_x} \mbox{ has a nonzero coefficient in } p(t,x)\right\}.
	\end{equation*}
	we may visualize it via the \emph{Newton polytope} $\conv (\mathcal N(p))$, i.e. the convex polytope in $\R^2$ spanned by $\mathcal N(p)$ (see \cref{figNewtonP}). We can see that multiplication by a monomial corresponds to a translation of $\conv (\mathcal N(p))$. If we denote by $\pi_t: \R^2 \to \R$ the projection to the $t$-axis, then the length $l_t := \len (\pi_t(\conv (\mathcal{N}(p))))$ is invariant under translation and thus invariant under $\doteq$-equivalence, and clearly $l_t \leq \deg p^N$. We can define $\pi_x$ and $l_x$ analogously.
	
	Now consider the Alexander polynomial $\Delta_L^N(t,x)$. We assumed that $\deg_x(\Delta_L^N) \geq 1$, so $l_x \geq 1$ for the associated Newton polytope. As $\Delta_L^N$ is in normal form, the Newton polytope touches the $t$-axis, i.e. there is some $(a,0) \in \mathcal{N}(\Delta_L^N)$. As $l_x \geq 1$, there is also some $(b,c) \in \mathcal{N}(\Delta_L^N)$ with $c \geq 1$.
	
	We know by \cref{bicAlex} that $\Delta_{d^k(L)}^N(t,x) \doteq \Delta_L^N(t,xt^{wk})$, and also
	\begin{equation}
		\mathcal{N}\left( \Delta_L^N(t,xt^{wk}) \right) = \left\{ (m + nwk, n) \, \middle| \, (m,n) \in \mathcal{N}\left( \Delta_L^N(t,x) \right) \right\}.
	\end{equation}
	In particular we know that $(a,0), (b + cwk, c) \in \mathcal{N}(\Delta_L^N(t,xt^{wk}))$. Hence for the associated Newton polytope we have that $l_t \geq b + cwk - a$. But this polytope is a translation of the Newton polytope $\conv \left( \mathcal{N} \left(\Delta_{d^k(L)}^N\right)\right)$, thus we obtain that $\deg \Delta_{d^k(L)}^N \geq l_t \geq cwk + b - a$. In particular, since $w,c \geq 1$, this goes to infinity as $k \to \infty$, and there is some $k_+ \in \N$ such that $\deg \Delta_{d^k(L)}^N > \deg \Delta_L^N$ for all $k > k_+$.
	
	Analogously we observe that $l_t \geq a - (b + cwk)$ goes to $\infty$ as $k \to - \infty$, so there is also some $k_- \in \N$ such that $\deg \Delta_{d^k(L)}^N > \deg \Delta_L^N$ for all $k < k_-$.
	
	We have now found a finite set $[k_-,k_+] \subset \Z$ such that, if $k \notin [k_-, k_+]$, then $\deg \Delta_{d^k(L)}^{N+} > \deg \Delta_L^N(t,x)$ and thus $\Delta_{d^k(L)}^N > \Delta_L^N(t,x)$. The set
	\begin{equation}
		\left\{ \Delta_{d^k(L)}^N(t,x) \, \middle| \, k \in \left[k_-, k_+\right]\right\}
	\end{equation}
	is finite and thus has a minimal element. This minimal element must then also be a minimal element of $D\Delta(L)$.
\end{proof}
We choose this minimal element as our normalization under Dehn twists. Define the \emph{Dehn-reduced Alexander polynomial} as $\overline{\Delta}_L(t,x) := \min D\Delta(L)$. Note that the bounds $k_-, k_+$ in the proof depend only upon the polynomial $\Delta_L^N(t,x)$, so the reduced polynomial can be computed easily from the non-reduced one.
\begin{prop}
	\label{uniqueTwist}
	If $L$ has winding number $w > 0$ and for the Alexander polynomial, we have that $d_x = \deg_x(\Delta_L^N(t,x)) > 0$, then there is a unique $k \in \Z$ such that $\Delta_{d^k(L)}^N(t,x) = \overline{\Delta}_L(t,x)$.
\end{prop}
\begin{proof}
	The existence of such a $k$ is clear. Now assume we have some $k,l$ with $\Delta_{d^k(L)}^N(t,x) = \overline{\Delta}_L(t,x) = \Delta_{d^l(L)}^N(t,x)$. Then we have that $d^l(L) = d^{l-k}(d^k(L))$, so by \cref{bicAlex}
	\begin{equation}
		\label{eqMonoThmAppl}
		\Delta_{d^k(L)}^N(t,x) = \Delta_{d^l(L)}^N(t,x) = \pm t^q\Delta_{d^k(L)}^N(t,xt^{(l-k)w}).
	\end{equation}
	If we write
	\begin{equation}
		\Delta_{d^k(L)}^N(t,x) = \sum_{j=0}^{d_x} p_j(t)x^j,
	\end{equation}
	then \cref{eqMonoThmAppl} becomes
	\begin{equation}
		\sum_{j=0}^{d_x} p_j(t)x^j = \sum_{j=0}^{d_x} t^{(l-k)wj+q}p_j(t)x^j.
	\end{equation}
	Comparing the $\Lambda[t]$ coefficients $p_j(t)$, we thus obtain $p_0(t) = t^qp_0(t)$, i.e. $q=0$ (as before we know that $p_0(t) \neq 0$). Furthermore $p_{d_x}(t) \neq 0$, and there the comparison yields $p_{d_x}(t) = t^{(l-k)wd_x}p_{d_x}(t)$, and as $wd_x > 0$ this implies $k=l$, as desired.
\end{proof}

By construction, $\overline{\Delta}_L(x,t)$ is invariant under Dehn twists, and due to \cref{revAlex}, it is also invariant under $q$.

If we start with some unoriented link $L \in V^3$, we have $2^k$ different ways to put an orientation on $L$, where $k$ is the number of link components. We would like to reduce the number of orientations we need to consider. For an orientation $\mathcal{O}$ on $L$, let $L^\mathcal{O}$ denote the associated oriented link. In the following, $w$ denotes the total winding number of an oriented link.

\begin{definition}
	Let $L$ be an unoriented link in $V^3$.
	\begin{enumerate}
		\item The \emph{maximal winding number} of $L$ is
		\begin{equation}
			w_{\max}(L) = \max \left\{ w\left(L^\mathcal{O}\right) \, \mid \, \mathcal{O} \mbox{ is an orientation for } L \right\}.
		\end{equation}
		\item Let $\mathcal{O}$ be an orientation for $L$. We say that $\mathcal{O}$ is \emph{admissible} if $w\left( L^\mathcal{O} \right) = \pm w_{\max}(L)$.
	\end{enumerate}
\end{definition}

By \cref{windingWrapping}, we know that $w(L^\mathcal{O}) \leq \wrap$, and the wrapping number does not depend on orientation, so $w_{\max}(L) \leq \wrap$.

With this, we can define the invariant
\begin{equation}
	\mathcal{A}_L = \left\{ \overline{\Delta}_{L^\mathcal{O}}(t,x) \, \middle| \, \mathcal{O} \mbox{ admissible }\right\}.
\end{equation}

\begin{corollary}
	$\mathcal{A}_L$ is an invariant of h-equivalence.
\end{corollary}
\begin{proof}
	Since the Alexander polynomial is an a-equivalence invariant, it is clear that $\mathcal{A}_L$ is invariant under a-equivalence. Thus it suffices to prove that $\mathcal{A}_L = \mathcal{A}_{h(L)}$ for the generators $d,r,q \in \MCG^{\pm} (V^3)$. For $d$ and $q$, we have already seen that $\overline{\Delta}_{L^\mathcal{O}}(x,t)$ is invariant.\\
	It remains to show that $\mathcal{A}_L = \mathcal{A}_{r(L)}$. As before, note that the extension of the rotation $\overline r$ switches only the orientation of the augmentation component $m$. But \cref{alexSym} implies the Alexander polynomial is invariant under switching all orientations at once, so we may instead switch all orientations except the one on $m$, which would mean switching all orientations on $L$. In other words, we get
	\begin{equation*}
		\Delta_{L^{-\mathcal{O}}}(t,x) \doteq \Delta_{r(L^\mathcal{O})}(t,x),
	\end{equation*}
	where $-\mathcal{O}$ is the reserve orientation. Now, we have that $w(L^\mathcal{O}) = - w(L^{-\mathcal{O}})$ for the winding number, so if $\mathcal{O}$ was admissible, so is $-\mathcal{O}$. Hence, the normalization $\overline{\Delta}_{r(L^\mathcal{O})}(t,x) \in \mathcal{A}_L$.
\end{proof}

\section{The dichromatic Kauffman polynomial}
There is another polynomial invariant that we can utilize to analyze solid torus links. It was first defined by Hoste and Przytycki in \cite{HP89}, and we follow that definition (with the slight variation that we allow an empty link; this makes the invariant more suitable for computation).

Recall that $D_{\circ}$ denotes the standard diagram for the affine unknot, and $D_|$ the standard diagram for the longitude knot. Further, let $D_1 \sqcup D_2$ denote the diagram consisting of the diagrams $D_1$ and $D_2$ next to one another.
\begin{definition}
	Let $L \subset V^3$ be a link, and $D$ a periodic diagram for $L$. The \emph{dichromatic Kauffman bracket} $\left\langle D \right\rangle$ is the invariant of $D$ defined by the following properties:
	\begin{enumerate}
		\item $\left\langle D_{\circ} \right\rangle = (- A^{-2} - A^2)$
		\item $\left\langle D_| \right\rangle = x$
		\item $\big\langle \input{tikz/skein-start.tex} \big\rangle = A \big\langle \input{tikz/skein-plus.tex} \big\rangle + A^{-1} \big\langle \input{tikz/skein-minus.tex} \big\rangle$
		\item $\left\langle D_1 \sqcup D_2 \right\rangle = \left\langle D_1 \right\rangle \left\langle D_2 \right\rangle$
		\item $\langle \cdot \rangle$ is invariant under $\Theta_2^\pm$ moves.
	\end{enumerate}
\end{definition}
Here, condition 3 should be read as a local smoothing of a crossing while the rest of the diagram stays the same; this can be applied at any crossing of a diagram. This type of equation is often called a \emph{skein relation}.

The bracket is not invariant under the $\Omega_1$ move; indeed we have
\begin{equation}
	\label{omega1kauff}
	(-A)^{-3} \left \langle \input{tikz/move-omega-kauff-1.tex} \right\rangle = \left\langle \input{tikz/move-omega-kauff-2.tex} \right\rangle = (-A)^3 \left\langle \input{tikz/move-omega-kauff-3.tex} \right\rangle
\end{equation}
It may not be immediately clear that the dichromatic bracket is well-defined even on the diagrammatic level. The reason behind this is that it is really an invariant of the annulus diagram; the well-definedness of the annulus invariant was proven in \cite{HP89}.
\begin{lemma}
	The dichromatic Kauffman bracket is invariant under the moves $\Omega_2, \Omega_3$ and $\Theta_1^\pm$.
\end{lemma}
\begin{proof}
	The invariance under $\Omega_2, \Omega_3$ can be shown completely analogously to the invariance of the classical Kauffman bracket, see \cite[lemmas 2.3 and 2.4]{ka87}.
	
	$\Theta_1^\pm$ invariance follows by resolving the crossing via the Skein relation and using $\Theta_2^\pm$ moves. We show the $\Theta_1^+$ case; the other one is analogous.
	\begin{multline}
		\left\langle \input{tikz/move-theta-kauff-1.tex} \right\rangle = A \left\langle \input{tikz/move-theta-kauff-2.tex} \right\rangle + A^{-1} \left\langle \input{tikz/move-theta-kauff-3.tex} \right\rangle \\ \stackrel{(*)}{=} A \left\langle \input{tikz/move-theta-kauff-2.tex} \right\rangle + A^{-1} \left\langle \input{tikz/move-theta-kauff-4.tex} \right\rangle = \left\langle \input{tikz/move-theta-kauff-5.tex} \right\rangle
	\end{multline}
	Here $(\ast)$ follows via one $\Theta_2^+$ move followed by a $\Theta_2^-$ move.
\end{proof}
A useful tool for analyzing the Kauffman bracket are so-called \emph{Kauffman states}: Given a diagram $D$, a Kauffman state $\sigma$ is an assignment of one of the symbols $0$ or $\infty$ to each crossing. The resolution $D_\sigma$ of $D$ with respect to $\sigma$ is then obtained by replacing each crossing $\input{tikz/skein-start.tex}$ with $\input{tikz/skein-plus.tex}$ if the symbol is $0$ and with $\input{tikz/skein-minus.tex}$ if the symbol is $\infty$, resulting in a diagram with no crossings. The \emph{signature} $\sgn \sigma$ of a Kauffman state $\sigma$ is the number of $0$s minus the number of $\infty$s in $\sigma$. We may then write
\begin{equation}
	\label{eqKauffmanState}
	\langle D \rangle = \sum_{\sigma} A^{\sgn \sigma} \langle D_\sigma \rangle,
\end{equation}
and the Kauffman bracket of a diagram without crossings is $(-A^2 - A^{-2})^jx^k$, where $j$ is the number of affine components and $k$ is the number of non-affine components.

In order to turn the bracket into a true knot invariant, we need to normalize to account for the behavior under $\Omega_1$ moves. To this end, for an oriented periodic link diagram $D$ and a crossing $c$ in it, we say that the \emph{sign} of $c$ is $+1$ if the outgoing undercrossing strand lies to the left of the overcrossing strand, and $-1$ otherwise. The \emph{self-writhe} $\sw(D)$ is the sum over the signs of all crossings in which over- and undercrossing strand belong to the same link component. This does not actually depend on the orientation chosen: If we switch the orientation of one component, the signs of all crossings where it crosses itself remain the same. Thus it makes sense to speak of $\sw(D)$ for an unoriented diagram $D$.
\begin{definition}
	The \emph{dichromatic Kauffman polynomial} for a solid torus link $L$ with periodic diagram $D$ is defined as $\nabla_L (A, x) = (-A^3)^{\sw(D)} \left\langle D \right\rangle$.
\end{definition}
This is indeed well-defined since our normalization makes sure that $\nabla_L(A,x)$ is invariant under $\Omega_1$ moves; this is due to the fact that adding a kink changes the self-writhe by $+1$ or $-1$, and due to \cref{omega1kauff}. Thus the Kauffman polynomial gives an invariant of a-equivalence. For the rest of this section, we will study its behavior under h-equivalence.

\begin{prop}
	\label{kauffRot}
	For the rotation $r \in \MCG(V^3)$, we have that $\nabla_L(A,x) = \nabla_{r(L)}(A,x)$.
\end{prop}
\begin{proof}
	As seen before, $r$ corresponds to a $180^\circ$ rotation on the diagrammatic level; call this move $\widetilde r$. It is clear that the rotation does not change the self-writhe as it keeps the orientation of the crossings intact. Thus we only need to show that $\langle D \rangle = \langle \widetilde r (D) \rangle$ for any link diagram $D$. To see this, simply note that the rotation does not affect the properties of the bracket; that is, $\widetilde r$ fixes $D_\circ$ and $D_|$, and applying the Skein relation as in property 3 or the diagram split as in property 4 also commutes with the rotation $\widetilde r$, as does performing a $\Theta_2$ move.
\end{proof}
\begin{prop}
	For the reflection $q \in \MCG^{\pm}(V^3)$, we have that $\nabla_L(A, x) = \nabla_{q (L)} (A^{-1}, x)$.
\end{prop}
\begin{proof}
	We will instead prove that $\nabla_L(A,x) = \nabla_{q  \circ r (L)} (A^{-1},x)$; then the claim follows from \cref{kauffRot}. Since $r$ is a rotation of the round solid torus in $\R^3$, and $q$ a reflection across a plane incident with the rotation axis, the composition is a reflection across the $90^\circ$ rotation of that plane. It follows that on the annulus diagram (and thus also on the periodic diagram), the homeomorphism $q \circ r$ corresponds to flipping over all crossings from under to over.
	
	Let $D$ be a diagram of $L$ and $D^-$ the diagram with all crossings reversed. For a Kauffman state $\sigma$ of $D$, we may define the Kauffman state $\overline \sigma$ of $D^-$ by changing the symbols at all crossings. Then we get $D_\sigma = D^-_{\overline \sigma}$ and $\sgn \sigma = - \sgn \overline \sigma$. It follows that
	\begin{equation}
		\langle D \rangle = \sum_{\sigma} A^{\sgn \sigma} \langle D_\sigma \rangle = \sum_{\sigma} (A^{-1})^{\sgn \overline \sigma} \langle D^-_{\overline \sigma} \rangle,
	\end{equation}
	and since for diagrams without crossings we have
	\begin{equation}
		\left\langle D_{\sigma} \right\rangle (A, x) = (-A^2 - A^{-2})^jx^k = \left\langle D_{\sigma} \right\rangle (A^{-1}, x),
	\end{equation}
	we may conclude that $\left\langle D \right\rangle (A, x) = \left\langle q \circ r (D) \right\rangle (A^{-1}, x)$. Also, given any orientation, flipping a crossing will change its sign, which implies that the self-writhe of $D$ is minus the self-writhe of $D^-$. Altogether we obtain $\nabla_L(A,x) = \nabla_{q \circ r (L)} (A^{-1}, x)$.
\end{proof}
Thus we can see that, if we define an equivalence relation $p(A,x) \sim p(A^{-1},x)$, then the equivalence class $[\nabla_L(A,x)]$ is a symmetry invariant.
\begin{prop}
	\label{propWrappingParity}
	Let $D$ be a periodic diagram with wrapping number $\wrap$. If $\wrap$ is even, then $\langle D \rangle(A,x)$ has only even powers of $x$, and if $\wrap$ is odd, then it has only odd powers of $x$.
\end{prop}
\begin{proof}
	The Kauffman state sum \cref{eqKauffmanState} of $\langle D \rangle$ has only diagrams without crossings of wrapping number $\wrap$. After a sequence of $\Theta_2$ moves, each of these becomes a diagram with only circles and straight lines, and there the wrapping number is equal to the power of $x$. But each $\Theta_2$ move changes the wrapping number only by 2, so the parity doesn't change.
\end{proof}
Of course we may conclude immediately that the same statement holds for the Kauffman polynomial of a link with a given wrapping number.

For the next property, note that the \emph{leading term} of a polynomial $p \in \Z[x]$ is the summand with the highest $x$-power appearing in $p$.
\begin{prop}
	\label{propBracketNonzero}
	Let $D$ be a $k$-component periodic diagram with $j$ components of winding number zero. Then the leading term of $\langle D \rangle (1,x)$ is given by the formula $(-1)^{\sw(D) - w_{\max}(L) + k }2^jx^{w_{\max}(L)}$.
\end{prop}
\begin{proof}
	Consider the skein relations for the bracket if we switch a crossing:
	\begin{align}
		\label{eqSkeinPlus}
		\langle \input{tikz/skein-start.tex} \rangle & = A \langle \input{tikz/skein-plus.tex}\rangle + A^{-1} \langle \input{tikz/skein-minus.tex}\rangle, \\
		\label{eqSkeinMinus} \langle \input{tikz/skein-invert.tex} \rangle & = A \langle \input{tikz/skein-minus.tex}\rangle + A^{-1} \langle \input{tikz/skein-plus.tex}\rangle.
	\end{align}
	We may solve the \cref{eqSkeinMinus} for $\big\langle \input{tikz/skein-minus.tex} \big\rangle$ and plug this into \cref{eqSkeinPlus}; then after rearranging we get
	\begin{equation}
		(A - A^{-2}) \big\langle \input{tikz/skein-plus.tex} \big\rangle = \big\langle \input{tikz/skein-start.tex} \big\rangle - A^{-1} \big\langle \input{tikz/skein-invert.tex} \big\rangle.
	\end{equation}
	We may evaluate at $A = 1$; then by the above we have $\big\langle \input{tikz/skein-start.tex} \big\rangle (1,x) = \big\langle \input{tikz/skein-invert.tex} \big\rangle (1,x)$; i.e. a crossing change does not affect the value $\langle D \rangle (1,x)$. Thus $\langle D \rangle (1,x)$ is an invariant of the \emph{shadow} $\overline D$ of $D$; the diagram with crossing information removed. Write $\langle \overline D \rangle (x) := \langle D \rangle (1,x)$. As $\langle \cdot \rangle$ is invariant under $\Omega_2, \Omega_3, \Theta_1$ and $\Theta_2$ moves, so is $\langle \overline D \rangle$ invariant under the corresponding shadow moves $\overline \Omega_2, \overline \Omega_3, \overline \Theta_1$ and $\overline \Theta_2$. The move $\Omega_1$ changes the Kauffman bracket by a factor of $(-A)^{\pm 3}$ and thus changes the sign of $\langle \overline D \rangle$.
	
	Any shadow $\overline D$ may be reduced to the shadow of a basic diagram: If two segments in $\overline D$ intersect twice, we can perform a sequence of $\overline \Omega_3$ and $\overline \Omega_2$ moves to eliminate that double crossing, and for self-crossings we may need an $\overline \Omega_1$ move. We may reduce further: Any segment that starts and ends on the same boundary component can then be removed with a sequence of $\overline \Omega_3$ and $\overline \Theta_1$ moves followed by a final $\overline \Theta_2$ move.
	
	This means that $\overline D$ is equivalent to a shadow $\overline{D'}$ of some basic diagram $D'$ consisting only of closed segments and segments starting and ending on opposite sides. $\overline{D'}$ is completely determined by a permutation $\sigma$ on the passings, and the number $c$ of closed segments. We get one closed segment for each component of $L$ with winding number $0$, and the number of passings of $\overline{D'}$ is equal to the maximal winding number $w$ of $L$. Let $s$ denote the number of $\overline \Omega_1$ moves used. Then $\langle \overline D \rangle = (-1)^s \langle \overline{D'} \rangle$.
	
	Any shadow move other than $\overline \Omega_1$ changes the number of crossings by a multiple of 2 if at all. Thus, the difference in the number of crossings between $D$ and $D'$ is odd if and only if $s$ is odd. In a basic diagram determined by $\sigma$ and $c$, any component with a winding number $w_0 > 0$ corresponds to a cycle of length $w_0$ in $\sigma$; such a cycle can be written as the product of transpositions, and $\overline{D'}$ corresponds to such a factorization: each self-crossing gives a transposition.
	
	A cycle of length $w_0$ has signature $(-1)^{w_0 -1}$ and thus the associated component in $\overline{D'}$ has $w_0 - 1 + 2n$ self-crossings for some $n \in \Z$. Component winding numbers do not change under shadow moves, and so $D'$ has $k - j$ components with winding number $> 0$ and thus the number of crossings in $D'$ is $w_{\max}(L) - (k-j) + 2m$ for some $m \in \Z$.
	
	So if $c$ is the number of crossings in $D$, then $c - w_{\max}(L) + k - j$ has the same parity as $s$. Crossings between different components always appear in pairs, so we may replace $c$ by only the self-crossings, and the number of self-crossings in turn has the same parity as $\sw D$; overall we can conclude that
	\begin{equation}
		\label{eqParityBracket}
		\langle \overline D \rangle = (-1)^{\sw D - w_{\max}(L) + k - j} \langle \overline{D'}\rangle.
	\end{equation}
	
	We note that there is exactly one Kauffman state $\sigma_{\max}$ of $D'$ such that the maximal winding number of $D'_{\sigma_{\max}}$ is equal to that of $D'$: the one where all crossings are resolved in such a way that no segment is created that starts and ends on the same boundary component.
	
	Then $\langle D'_{\sigma_{\max}} \rangle = (-A^{2}-A^{-2})^j x^{w_{\max}(L)}$, and for any other Kauffman state $\sigma'$ only powers of $x$ smaller than $w_{\max}(L)$ appear in $\langle D'_{\sigma'} \rangle$; since $\langle D' \rangle$ is the sum of all Kauffman states, we may conclude that the leading term of $\langle \overline D' \rangle$ is $(-2)^jx^{w_{\max}(L)}$. From this and \cref{eqParityBracket} we may conclude that the leading term of $\langle \overline D \rangle$ is $(-1)^{\sw D - w_{\max}(L) + k}2^jx^{w_{\max}(L)}$.
\end{proof}
\begin{corollary}
	\label{corollKauffNonzero}
	Let $L$ be a $k$-component solid torus link with $j$ components of winding number zero. Then the leading term of $\nabla_L(1,x)$ is given by the formula $(-1)^{k-w_{\max}(L)} 2^jx^{w_{\max}(L)}$. In particular we have that $\nabla_L(A,x) \neq 0$.
\end{corollary}
\begin{proof}
	We have that $\nabla_L(A,x) = (-A^3)^{\sw D} \langle D \rangle$ for any periodic diagram $D$ of $L$. Thus $\nabla_L(1,x) = (-1)^{\sw D} \langle D \rangle$, and the statement follows via \cref{propBracketNonzero}.
\end{proof}

To conclude our general remarks on the structure of the Kauffman polynomial, we note a further simple property:
\begin{prop}
	\label{propWrappingBound}
	Let $L$ be any solid torus link with wrapping number $\wrap$. Then $\deg_x(\nabla_L) \leq \wrap$.
\end{prop}
\begin{proof}
	If $D$ is a diagram of $L$ with $\wrap$ passings, then for any Kauffman state $\sigma$ of $D$, $D_\sigma$ also has $\wrap$ passings, and thus can have at most $\wrap$ longitude knot components. The statement follows via \cref{eqKauffmanState}.
\end{proof}
Indeed equality holds for all links in our tabulation. We suspect it may hold in general.
\begin{conjecture}
	\label{conjWrapping}
	Let $L$ be any solid torus link with wrapping number $\wrap$. Then $\deg_x(\nabla_L) = \wrap$.
\end{conjecture}
We can indeed show that this is true whenever $\wrap \leq 1$: By \cref{corollKauffNonzero}, we know that $\deg_x(\nabla L) \geq 0$, so there is some nontrivial coeffient, and then the statement follows due to \cref{propWrappingParity}. This however does not generalize to larger wrapping numbers.

\begin{wrapfigure}{r}{0.2\textwidth}
	\begin{center}
		\input{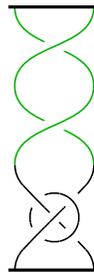}
		\caption{The product $D$\textcolor{black!30!green}{$\otimes D_\Delta$} for a wrapping number 2 diagram.}
		\label{dehntwist}
	\end{center}
\end{wrapfigure}
Next, we study the effect of a Dehn twist on the Kauffman polynomial. Here our main result will be that for small wrapping numbers, the dichromatic Kauffman polynomials of two links provide an upper bound for the degree that a Dehn twist taking one link into the other can have.

First of all, let us note a few things: It is clear that for links with wrapping number 0 (i.e. affine links) or 1, the Dehn twist has no effect on the link; that is, in this case a link and its Dehn twist will be a-equivalent. The first case where we obtain interesting behavior is wrapping number 2.

\begin{figure}
	\input{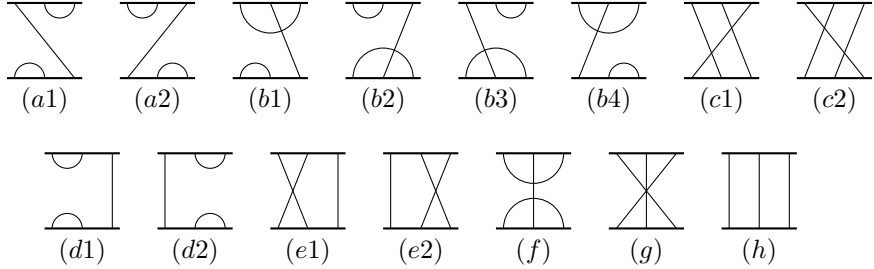}
	\caption{All configurations of endpoints for a diagram with wrapping number 3. They are sorted by the number of link components they have.}
	\label{figEndOptions}
\end{figure}

For the remainder of this chapter, it will be useful to define an operation on diagrams: Let $D_1, D_2$ be two periodic diagrams of equal wrapping number; then we can construct a new diagram $D_1 \otimes D_2$ by stacking $D_2$ on top of $D_1$, gluing the upper passings of $D_1$ to the lower passings of $D_2$. Note that $D_1 \otimes D_2$ can be transformed into $D_2 \otimes D_1$ via a series of $\Theta_i$ moves; in that sense our operation is almost commutative. We will also adopt the notation
\begin{equation}
	D^k = \underbrace{D \otimes \ldots \otimes D}_{k \mbox{ times}}.
\end{equation}
Let $D$ be a periodic diagram of a link $L$ with wrapping number 2, and let $\overline{D}_k$ be the $k$-fold Dehn twist of the diagram $D$ (i.e. the $k$-fold application of the move $\Delta$ as in \cref{hequivmoves}) for some $k \in \N$. Then we can write $\overline D_k = D \otimes D_\Delta^k$, where $D_\Delta$ is a diagram of the twisted region (with 2 strands).
\begin{lemma}
	\label{kauffdehn2}
	For the Kauffman bracket $\langle \cdot \rangle$ we have the formula
	\begin{equation}
	\left\langle \overline{D}_k \right\rangle = \left( 1 - A^{-4} \right) \sum_{i=0}^{k-1} A^{2(k-1) - 8i} \left\langle D\otimes \left( \cupdiag \right) \right\rangle + A^{2k} \left\langle D \right\rangle
	\end{equation}
	if $k \geq 0$.
\end{lemma}
\begin{proof}
	The formula trivially holds for $k = 0$. For $k=1$ we can compute:
	\begin{equation}
	\label{kdehn2indstart}
	\begin{aligned}
	\left\langle \overline{D}_{1} \right\rangle = & \left\langle D \otimes D_\Delta \right\rangle = \left\langle D \otimes \left( \input{tikz/kauff-dehn-2-proof-1.tex} \right)\right\rangle \\ = & A \left\langle D \otimes \left( \input{tikz/kauff-dehn-2-proof-2.tex} \right) \right\rangle + A^{-1} \left\langle D \otimes \left( \input{tikz/kauff-dehn-2-proof-3.tex} \right) \right\rangle \\
	= & A^2 \left\langle D \right\rangle + \left(1 - A^{-4} \right) \left\langle D \otimes \left( \cupdiag \right) \right\rangle
	\end{aligned}
	\end{equation}
	The final equality is due to the skein relation and \cref{omega1kauff}. Now by induction we may assume the formula holds for some $k > 0$. Then
	\begin{equation}
		\begin{aligned}
			\left\langle \overline{D}_{k+1} \right\rangle = & \left\langle \overline{D}_k \otimes D_\Delta \right\rangle \stackrel{\cref{kdehn2indstart}}{=} A^2 \left\langle \overline{D}_k \right\rangle + (1-A^{-4}) \left\langle \overline{D}_k \otimes \left( \cupdiag \right) \right\rangle \\
			\stackrel{(*)}{=} & A^2 \left( \left( 1 - A^{-4} \right) \sum_{i=0}^{k-1} A^{2(k-1) - 8i} \left\langle D \otimes \left( \cupdiag \right) \right\rangle + A^{2k} \left\langle D \right\rangle \right) \\
			& + \left( 1 - A^{-4}\right) \left\langle D \otimes D_\Delta^k \otimes \left( \cupdiag \right) \right\rangle \\
			\stackrel{(**)}{=} & \left( 1 - A^{-4} \right) \left( \sum_{i=0}^{k-1} A^{2k - 8i} \left\langle D \otimes \left( \cupdiag \right) \right\rangle + A^{2(k+1)} \left\langle D \right\rangle + A^{-6k} \left\langle D \otimes \left( \cupdiag \right) \right\rangle \right) \\
			= & \left( 1 - A^{-4} \right) \sum_{i=0}^{k} A^{2k - 8i} \left\langle D \otimes \left( \cupdiag \right) \right\rangle + A^{2(k+1)} \left\langle D \right\rangle
		\end{aligned}
	\end{equation}
	Here the equality $(*)$ holds due to the induction assumption and $(**)$ holds via \cref{omega1kauff} since $D \otimes D_\Delta^k \otimes \left( \cupdiag \right)$ is related to $D \otimes \left( \cupdiag \right)$ by a series of $2k$ Reidemeister 1 moves.
\end{proof}
For what follows, we need a notion of degree of the Kauffman polynomial; for this note that by definition, $\nabla_L(A,x)$ cannot have any negative powers in $x$; thus the $x$-degree $\deg_x(\nabla_L)$ is well-defined. On the other hand, $A$ may have negative powers, so we need to adjust our definition.
\begin{definition}
	Let $p(x_1, \ldots, x_k)$ be a nonzero Laurent polynomial in $k$ variables. The \emph{$x_i$-spread} $\spread_{x_i} (p)$ is the difference between the maximal and the minimal power of $A$ occurring in $p$. The \emph{spread} $\spread (L)$ of a solid torus link $L$ is the $A$-spread of its Kauffman polynomial.
\end{definition}
The $x_i$-spread is well-behaved under multiplication, we have:
\begin{equation}
	\spread_{x_i} (pq) = \spread_{x_i} (p) + \spread_{x_i} (q)
\end{equation}
It is also obvious that the spread of a link $L$ equals the $A$-spread of the Kauffman bracket of any diagram of $L$.
\begin{definition}
	The \emph{$x$-capped spread} $\widehat{\spread}(p)$ of a nonzero Laurent polynomial $p(A,x)$ is defined as the difference between the smallest power of $A$ occuring and the largest power of $A$ occuring in a monomial wherein the power of $x$ is also maximal.
	
	The $x$-capped spread of a solid torus link is defined as the $x$-capped spread of its Kauffman polynomial.
\end{definition}
Note that by \cref{corollKauffNonzero}, $\nabla_L$ is always nonzero, so the $x$-capped spread $\widehat{\spread}(L)$ is always defined; there is always some smallest power of $A$.
\begin{lemma}
	\label{lemmaWrap2}
	Let $D$ be a periodic diagram of a link $L$ with wrapping number 2 such that $\deg_x(\nabla_L) = 2$, and further let $a_0$ be the smallest power of $A$ appearing in the bracket polynomial $\langle D \rangle (A,x)$ and $a_1$ be the smallest power of $A$ appearing in $\langle D \otimes \left( \cupdiag \right) \rangle$. Assuming $L'$ is isotopic to a $k$-fold Dehn twist of $L$ for some $k \geq 0$, we have the following bound for $k$:
	\begin{equation*}
		8k \leq a_1 + 2 - a_0 + \max \left(0, \widehat \spread\left(L'\right) - \widehat \spread(L) \right)
	\end{equation*}
\end{lemma}
\begin{proof}
	We assume that $8k > a_1 + 2 - a_0$, and wish to show that then, $8k \leq \widehat \spread\left(L'\right) - \widehat \spread(L) + a_1 + 2 - a_0$ (indeed we will show equality in this case).
	
	Consider that $\widehat \spread(L') = \widehat \spread(\langle \overline{D}_k\rangle)$, and using \cref{kauffdehn2} we may conclude that
	\begin{equation}
		\widehat \spread(L') = \widehat \spread \left( \left( 1 - A^{-4} \right) \sum_{i=0}^{k-1} A^{2(k-1) - 8i} \left\langle D\otimes \left( \cupdiag \right) \right\rangle + A^{2k} \left\langle D \right\rangle \right).
	\end{equation}
	We note that $D \otimes \left( \cupdiag \right)$ is a diagram of an affine link; it follows that $\deg_x\left(\left\langle D \otimes \left( \cupdiag \right) \right\rangle \right) = 0$. On the other hand, by assumption $\deg_x(\langle D \rangle) = 2$; this implies that the largest $x$ power of $\langle \overline{D}_k\rangle$ needs to appear in $A^{2k}\langle D \rangle$, and the largest power of $A$ occuring in a monomial which contains $x^2$ is exactly $2k + \widehat \spread(L) + a_0$.
	
	Consider the lowest $A$-power appearing in the RHS. It is either the lowest power of the second summand, i.e. $2k + a_0$, or the lowest $A$-power of the first summand - unless these powers are equal, in which case they may cancel. Both summands are nonzero as the Kauffman bracket is never zero. The lowest power in the first summand is easily seen to be $-6k+2+a_1$. As we assumed that $8k > a_1 + 2 -a_0$, it follows that that power is strictly lower than the one in the second summand.
	
	Overall we conclude that
	\begin{equation}
		\widehat \spread(L') = 2k + \widehat \spread \langle D \rangle + a_0 - (-6k + 2 + a_1) = 8k + \widehat \spread \langle D \rangle - 2 + a_0 - a_1,
	\end{equation}
	and the claim follows.
\end{proof}
This bound unfortunately only works for links with wrapping number 2. It is possible, but a bit more complicated, to obtain a similar result for link with wrapping number 3. For this, it is necessary to first introduce a technical lemma. Towards this, note that we may classify diagrams of wrapping number 3 according to which passings are connected by a segment; there are 15 distinct possibilities as outlined in \cref{figEndOptions}, the number of matchings of 6 elements.
\begin{lemma}
	Let $D$ be a periodic diagram of any type other than $(e1)$, $(e2)$ or $(g)$, then we have that
	\begin{multline}
		\left( A^4 - A^{-4}\right) \left( \left\langle D \otimes \left( | \cupdiag \right) \right\rangle + \left\langle D \otimes \left( \cupdiag | \right) \right\rangle \right) \\ + \left( A^2 - A^{-2} \right) \left( \left\langle D \otimes \left( \input{tikz/kauff-dehn-3-basic-1.tex} \right) \right\rangle + \left\langle D \otimes \left( \input{tikz/kauff-dehn-3-basic-2.tex} \right) \right\rangle \right) \neq 0.
	\end{multline}
\end{lemma}
\begin{proof}
	As shorthand, we write
	\begin{alignat*}{2}
		D_{10} & := D \otimes \left( | \cupdiag \right), \quad & D_{01} & := D \otimes \left( \cupdiag | \right), \\
		D_\% & := D \otimes \left( \input{tikz/kauff-dehn-3-basic-1.tex} \right), \quad & D_\revPer & := D \otimes \left( \input{tikz/kauff-dehn-3-basic-2.tex} \right).
	\end{alignat*}
	Assume towards a contradiction that the sum equals zero for some periodic diagram $D$. As $(A^4 - A^{-4}) = (A^2 + A^{-2})(A^2 - A^{-2})$, we may divide the entire expression by $(A^2 - A^{-2})$ to obtain
	\begin{equation}
		\label{eqSigmaZero}
		\Sigma_D(A,x) := \left( A^2 + A^{-2}\right) \left( \left\langle D_{10} \right\rangle + \left\langle D_{01} \right\rangle \right) + \left\langle D_\% \right\rangle + \left\langle D_\revPer \right\rangle = 0.
	\end{equation}

	We consider evaluation at $A = 1$. Note that each of the four diagrams $D_*$ above admits a $\Theta_2$ move to a diagram with wrapping number 1; then by \cref{propWrappingBound} and \cref{propWrappingParity}, the only power of $x$ that can occur in $\langle D_* \rangle$ is $1$. By \cref{propBracketNonzero}, we can conclude $\langle D_* \rangle (1,x) = (-1)^{\sw D_*} (-2)^{j_*}$, where $j_*$ is the number of components of $D_*$ with winding number $0$.  In particular it follows that 
	\begin{multline}
		\label{eqSigmaVal1}
		\Sigma_D (1,x) = -2 \left( (-1)^{\sw D_{10}} (-2)^{j_{10}}x + (-1)^{\sw  D_{01}} (-2)^{j_{01}}x \right) \\ - (-1)^{\sw D_\%} (-2)^{j_\%}x - (-1)^{\sw D_\revPer} (-2)^{j_\revPer}x.
	\end{multline}
	We will first show that all of the self-writhes involved are the same mod 2, meaning that we can divide by $(-1)^{\sw D_*}$ to simplify the above expression. Thus we show that, for $a,b \in \{10, 01, \%, \revPer \}$, we have $\sw D_a - \sw D_b \equiv 0 \mod 2$.
	
	To this end, we count the signed self-crossings in the diagrams. Note that any self-crossing of a single segment in $D$ will still be a self-crossing of the same parity in all the $D_*$, so they will have no influence on the difference. We only need to consider the self-crossings of link components in the link associated to $D_*$ that involve different segments in $D$; these will be among the three segments that meet the boundary.
	
	For a pair of these segments, the number of their pairwise crossings mod 2 is determined by the relative positioning of their boundary points: there has to be an odd number if one segment has to cross the other to connect its boundary points, and an even number otherwise. We may use this observation to restrict ourselves to the 15 cases in \cref{figEndOptions}.
	
	We only need to consider the diagram $D$ up to vertical and horizontal mirror symmetry, as the expression for $\Sigma_D(A,x)$ is horizontally symmetric, and if $\widehat D$ is the vertical mirror image of $D$, then $D \otimes (| \cupdiag)$ can be transformed into $(| \cupdiag) \otimes D$, which is the vertical mirror image of $\widehat D \otimes (| \cupdiag)$ (and analogously for the other products). Up to such symmetry, there are only 8 distinct cases: $(a1)$, $(b1)$, $(c1)$, $(d1)$, $(e)$, $(f)$, $(g)$ and $(h)$. Of those, in the configurations $(a1)$, $(d1)$ and $(h)$ any two of the segments can only have an even number of crossings between them and thus the self-writhe mod 2 remains unaffected. We merely need to investigate the five remaining cases.
	 
	 For example, consider case $(c1)$. Label the segments $s_1, s_2, s_3$ from the bottom right. Then $\{ s_1, s_2\}$ and $\{ s_1, s_3\}$ have an odd number of crossings each, whereas $\{ s_2, s_3\}$ have an even number. In $D_{10}, D_{01}$ and $D_\%$, all the crossings are self-crossings, so we have an even number. In $D_\revPer$,  only the crossings of $\{ s_2, s_3\}$ are self-crossings, so there we have an even number as well.
	 
	 The other cases may be argued analogously in a completely straightforward fashion.
	 
	 Thus, from \cref{eqSigmaZero} and \cref{eqSigmaVal1} we conclude
	 \begin{equation}
	 	\pm \Sigma_D(1,x) = \left(2 \left( (-2)^{j_{D_{10}}} + (-2)^{j_{D_{01}}} \right) + (-2)^{j_{D_\%}} +  (-2)^{j_{D_\revPer}}\right)x = 0.
	 \end{equation}
 	We next investigate the terms $j_*$; for this we determine the number of loops. Note that any closed segment of $D$ becomes a loop in $D_*$; let the number of closed segments be $n_0$. To determine the additional loops we can again make a case distinction as in \cref{figEndOptions}. It is then straightforward to compute the $j_*$ and from that, the value $\Sigma_D(1,x)$. We again do this only up to symmetry.
 	\begin{center}
 	\begin{tabular}{c|c|c|c|c|c|c|c|c}
 		& \begin{tikzpicture}[scale=0.75]
 			\draw [thick] (0, 0) -- (1, 0);
 			\draw [thick] (0, 1) -- (1, 1);
 			\draw (0.5,1) arc (180:360:0.2);
 			\draw (0.5,0) arc (0:180:0.2);
 			\draw (0.9,0) -- (0.1,1);
 		\end{tikzpicture} & \begin{tikzpicture}[scale=0.75]
	 		\draw [thick] (0, 0) -- (1, 0);
	 		\draw [thick] (0, 1) -- (1, 1);
	 		\draw (0.1,1) arc (180:360:0.4);
	 		\draw (0.5,0) arc (0:180:0.2);
	 		\draw (0.9,0) -- (0.5,1);
	 	\end{tikzpicture} & \begin{tikzpicture}[scale=0.75]
	 		\draw [thick] (0, 0) -- (1, 0);
	 		\draw [thick] (0, 1) -- (1, 1);
	 		\draw (0.1,0) -- (0.9,1);
	 		\draw (0.5,0) -- (0.1,1);
	 		\draw (0.9,0) -- (0.5,1);
	 	\end{tikzpicture} & \begin{tikzpicture}[scale=0.75]
 			\draw [thick] (0, 0) -- (1, 0);
 			\draw [thick] (0, 1) -- (1, 1);
 			\draw (0.1,1) arc (180:360:0.2);
 			\draw (0.5,0) arc (0:180:0.2);
 			\draw (0.9,0) -- (0.9,1);
 		\end{tikzpicture} & \begin{tikzpicture}[scale=0.75]
	 		\draw [thick] (0, 0) -- (1, 0);
	 		\draw [thick] (0, 1) -- (1, 1);
	 		\draw (0.1,0) -- (0.5,1);
	 		\draw (0.5,0) -- (0.1,1);
	 		\draw (0.9,0) -- (0.9,1);
	 	\end{tikzpicture} & \begin{tikzpicture}[scale=0.75]
			 \draw [thick] (0, 0) -- (1, 0);
			 \draw [thick] (0, 1) -- (1, 1);
			 \draw (0.1,1) arc (180:360:0.4);
			 \draw (0.9,0) arc (0:180:0.4);
			 \draw (0.5,0) -- (0.5,1);
		\end{tikzpicture} & \begin{tikzpicture}[scale=0.75]
			\draw [thick] (0, 0) -- (1, 0);
			\draw [thick] (0, 1) -- (1, 1);
			\draw (0.1,0) -- (0.9,1);
			\draw (0.5,0) -- (0.5,1);
			\draw (0.9,0) -- (0.1,1);
		\end{tikzpicture} & \begin{tikzpicture}[scale=0.75]
			\draw [thick] (0, 0) -- (1, 0);
			\draw [thick] (0, 1) -- (1, 1);
			\draw (0.1,0) -- (0.1,1);
			\draw (0.5,0) -- (0.5,1);
			\draw (0.9,0) -- (0.9,1);
		\end{tikzpicture} \\
		$j_{01}$ & 1 & 1 & 0 & 2 & 1 & 0 & 0 & 1 \\
		$j_{10}$ & 1 & 0 & 0 & 0 & 0 & 0 & 0 & 1 \\
		$j_\%$ & 0 & 0 & 1 & 1 & 0 & 0 & 1 & 0 \\
		$j_{\revPer}$ & 2 & 1 & 0 & 1 & 0 & 0 & 1 & 0 \\
		$\displaystyle \pm\frac{\Sigma_D(1,x)}{(-2)^{n_0}x}$ & $-3$ & $-3$ & $3$ & $6$ & $0$ & $6$ & $0$ & $-6$
 	\end{tabular}
 	\end{center}
 	It thus follows that if $\Sigma_D(1,x) = 0$, the diagram $D$ must be of type $(e1)$, $(e2)$ or $(g)$, as desired.
\end{proof}
\begin{lemma}
	\label{kauffdehn3}
	Let $D$ be a periodic diagram of a link $L$ with wrapping number 3, and $\overline{D}_k$ the $k$-fold Dehn twist of $D$ for some $k \in \Z$. Then, if $k \geq 0$, we have
	\begin{multline}
		\label{eqKauffWrap3}
	\left\langle \overline{D}_k \right\rangle = \left( \left( A^4 - A^{-4}\right) \left( \left\langle D \otimes \left( | \cupdiag \right) \right\rangle + \left\langle D \otimes \left( \cupdiag | \right) \right\rangle \right) \right. \\ \left. + \left( A^2 - A^{-2} \right) \left( \left\langle D \otimes \left( \input{tikz/kauff-dehn-3-basic-1.tex} \right) \right\rangle + \left\langle D \otimes \left( \input{tikz/kauff-dehn-3-basic-2.tex} \right) \right\rangle \right) \right) \sum_{i=0}^{k-1} A^{6(k-1) - 12i} + A^{6k} \langle D \rangle.
	\end{multline}
\end{lemma}
\begin{proof}
	We will proceed as in the proof of \cref{kauffdehn2}. Let us first compute the case $k=1$.
	\begin{equation}
		\begin{aligned}
			\left\langle D \otimes D_\Delta \right\rangle = & \left\langle D \otimes \left( \input{tikz/kauff-dehn-3-proof-1.tex} \right) \right\rangle = A \left\langle D \otimes \left( \input{tikz/kauff-dehn-3-proof-2.tex} \right) \right\rangle + A^{-1} \left\langle D \otimes \left( \input{tikz/kauff-dehn-3-proof-3.tex} \right) \right\rangle \\
			 = & A^2 \left\langle D \otimes \left( \input{tikz/kauff-dehn-3-proof-4.tex} \right) \right\rangle + \left\langle D \otimes \left( \input{tikz/kauff-dehn-3-proof-5.tex} \right) \right\rangle + A^{-1} \left\langle D \otimes \left( \input{tikz/kauff-dehn-3-proof-6.tex} \right) \right\rangle \\
			 = & A^3 \left\langle D \otimes \left( \input{tikz/kauff-dehn-3-proof-7.tex} \right) \right\rangle + A \left\langle D \otimes \left( \input{tikz/kauff-dehn-3-proof-8.tex} \right) \right\rangle + \left\langle D \otimes \left( \input{tikz/kauff-dehn-3-proof-9.tex} \right) \right\rangle \\
			 & - A^{-4} \left\langle D \otimes \left( | \cupdiag \right) \right\rangle \\
			 = & A^4 \left\langle D \otimes \left( \input{tikz/kauff-dehn-3-proof-10.tex} \right) \right\rangle + A^2 \left\langle D \otimes \left( \input{tikz/kauff-dehn-3-proof-11.tex} \right) \right\rangle - A^{-2} \left\langle D \otimes \left( \input{tikz/kauff-dehn-3-basic-2.tex} \right) \right\rangle \\ & - A^{-3} \left\langle D \otimes \left( \input{tikz/kauff-dehn-3-proof-12.tex} \right) \right\rangle - A^{-4} \left\langle D \otimes \left( | \cupdiag \right) \right\rangle \\
			= & A^5 \left\langle D \otimes \left( \input{tikz/kauff-dehn-3-proof-13.tex} \right) \right\rangle + A^3 \left\langle D \otimes \left( \input{tikz/kauff-dehn-3-proof-14.tex} \right) \right\rangle + \left(A^2 - A^{-2}\right) \left\langle D \otimes \left( \input{tikz/kauff-dehn-3-basic-1.tex} \right) \right\rangle \\ & - A^{-2} \left\langle D \otimes \left( \input{tikz/kauff-dehn-3-basic-2.tex} \right) \right\rangle - A^{-4} \left( \left\langle D \otimes \left( \cupdiag | \right) \right\rangle + \left\langle D \otimes \left( | \cupdiag \right) \right\rangle \right) \\
			= & \left(A^4 - A^{-4}\right) \left( \left\langle D \otimes \left( \cupdiag | \right) \right\rangle + \left\langle D \otimes \left( | \cupdiag \right) \right\rangle \right)\\ & + \left( A^2 - A^{-2} \right) \left( \left\langle D \otimes \left( \input{tikz/kauff-dehn-3-basic-1.tex} \right) \right\rangle + \left\langle D \otimes \left( \input{tikz/kauff-dehn-3-basic-2.tex} \right) \right\rangle \right) + A^6 \left\langle D \right\rangle
		\end{aligned}
	\end{equation}
	Here we use the skein relation as well as $\Omega_i$ moves. The general claim follows analogously to \cref{kauffdehn2} since for each diagram $S \in \left\{ \cupdiag |, | \cupdiag, \input{tikz/kauff-dehn-3-basic-1.tex}, \input{tikz/kauff-dehn-3-basic-2.tex} \right\}$, we have that $D_\Delta \otimes S$ is related to $S$ via a sequence of two $\Omega_1$ as well as some number of $\Omega_2$ and $\Omega_3$ moves.
\end{proof}

\begin{lemma}
	\label{lemmaWrap3}
	Let $D$ be a periodic diagram of type other than $(e1)$, $(e2)$ or $(g)$ for a link $L$ with wrapping number 3 such that $\deg_x(\nabla_L) = 3$, and further let $a_0$ be the smallest power of $A$ appearing in the bracket polynomial $\langle D \rangle (A,x)$ and $a_1$ be the smallest power of $A$ appearing in $\Sigma_D(A,x)$ as on \cref{eqSigmaZero}. Further let $L'$ be isotopic to a $k$-fold Dehn twist of $L$ for some $k \geq 0$, then we have the following bound for $k$:
	\begin{equation*}
		12k \leq a_1 + 4 - a_0 + \max \left(0, \widehat \spread\left(L'\right) - \widehat \spread(L) \right)
	\end{equation*}
\end{lemma}
\begin{proof}
	Note that all diagrams appearing in $\Sigma_D$ are of links with wrapping number 1; then it follows that $\deg_x(\Sigma_D) \leq 1 < \deg_x( \langle D \rangle)$. With this, the proof proceeds completely analogously to that of \cref{lemmaWrap2}.
\end{proof}

\begin{theorem}
	Let $L$ and $L'$ be two links of equal wrapping number $\wrap \in \{ 2,3\}$ such that $\deg_x(\nabla_L) = \deg_x(\nabla_{L'}) = \wrap$, and $D, D'$ diagrams of $L, L'$ that realize the wrapping number. Then there is a $k_0 \in \N$ computable from the diagrams $D, D'$ such that if $L'$ is $h^+$-equivalent to $L$, then $L'$ has to be a $k$-fold Dehn twist of $L$ for some $| k | \leq k_0$.
\end{theorem}
\begin{proof}
	We may assume that $k \geq 0$; otherwise we can switch the roles of $L$ and $L'$ to get a positive Dehn twist. If $\wrap=2$, the conditions of \cref{lemmaWrap2} are satisfied and thus we can set $k_0$ as that bound.
	
	If $\wrap = 3$ and $D$ is not of one of the types $(e1)$, $(e2)$ or $(g)$ from \cref{figEndOptions}, we may use the bound from \cref{lemmaWrap3}. Finally, if $D$ is of one of the types (8), (12) or (14), that means that $L = L_1 \cup L_2$ where $L_i$ is a link with wrapping number $i$. Since $L$ and $L'$ are $h^+$-equivalent, there is a sublink $L_2'$ of $L'$ that is $h^+$-equivalent to $L_2$ via the same homeomorphism. Then $L_2'$ also has wrapping number 2 and thus we have reduced our problem to the case $\wrap=2$, for which we already have a bound.
\end{proof}
The assumption on $\deg_x$ of the polynomials would be superfluous if we could prove \cref{conjWrapping}. Unfortunately, it seems that there is no simple generalization to cases $\wrap \geq 4$, so we can only utilize this in some special cases.

\section{Classification of solid torus links}
We now have all the ingredients available to build a classification algorithm. Our goal will be to classify all non-split solid torus links up to a given crossing number $n$, and up to h-equivalence. We have opted to not include (h-)split links for two reasons: First, we would never be able to call our list complete, as there are infinitely many split links even with crossing number 0. Secondly, even if we avoided this issue by going only up to some maximal wrapping number and only considered non-affine split components, this would still cause a lot of bloat in our list, for only a limited gain, as those split links can be easily constructed from the non-split link table.

On the other hand, we will include composite links. These are generally fewer in number, and there are some nontrivial results in deciding whether taking the composition of a link with an affine piece along one link component gives the same link as taking it along another component, and whether the composition with a mirror image yields a symmetric link. 
\subsection{Generating and ordering basic diagrams}
\label{secGenOrder}
We know from \cref{passBound} that any non-split link with at most $n$ crossings possesses a basic diagram with at most $2n$ passings. We will generate our links by generating these basic diagrams. 

For this, note that for a basic diagram, which segments cross is completely determined by the boundary points of the segments. If we label our potential passings by integers $1, \ldots, k$, for some $k \leq 2n$, we may consider the set of all partitions of
\begin{equation*}
\left\{ 1, 2, \ldots, k, -1, -2, \ldots , -k \right\}
\end{equation*}
into sets with two elements to obtain the segments in a basic diagram. We can remove partitions that would give too many crossings as well as those that result in split diagrams. 

Although it is determined which segments cross, the order of the crossings within a segment is not necessarily fixed; if 3 segments all cross each other, then there is a choice to be made. For any partition, we collect all possible choices of orderings.

Finally, we need to specify at each crossing which segment is the overcrossing. For this, note that our labeling of the passings induces a canonical orientation on each segment from the boundary point with smaller label to the one with a larger label (here the signed labels are considered). With this orientation, fixing a sign for the crossing determines which segment crosses over. Choosing all combinations of signs gives us a base set of possible diagrams.

The next step is to convert our data, which currently consists of tuples $(S, C, O)$, where $S$ contains pairs of boundary points, $C$ orderings of the crossings along each segment, and $O$ an orientation for every crossing, into a format that can also handle non-basic diagrams. This is done using \emph{periodic Gauss paragraphs} (PGPs). A PGP of an $n$-crossing periodic diagram with $k$ passings consists of the following data:
\begin{itemize}
	\item An $n$-tuple $\omega$ of elements of $\{ -1, 1\}$, giving the orientations of the crossings,
	\item For each segment, a list $C_s$ of the form $(c_1, \ldots, c_m)$, where $c_i \in \{\pm 1, \ldots, \pm n\}$ indicate the crossings, with over-/undercrossing information,
	\item For each segment, a tuple $B_s$ that is either empty (if the segment has no boundary) or has two elements $(b_1, b_2)$ with $b_i \in \{ 1, \ldots, k, -1, \ldots, -k\}$ indicating the boundary points.
\end{itemize}
This is a generalization of an old notion originally due to Gauss, as a way of describing immersed planar curves (and by extension, diagrams of links in $S^3$). A more recent analysis of these (classical) Gauss paragraphs or codes focusing on computational efficiency can be found in \cite{RT84}.

We can directly convert our basic diagrams into PGPs if we give some labeling $1,\ldots, l$ to the crossings. Then, we use $\Theta_i$ moves to reduce the number of passings as much as possible; for now, we will keep the structure of the underlying graph of the annulus diagram fixed, thus we will not allow $\Omega_i$ moves at this stage.

As each link can have a variety of different diagrams and even the same diagram can be generated by different PGPs, we need a way to find a standard representation of a link. To this end, we introduce a total ordering on PGPs. This order is determined by checking two PGPs lexicographically for the following conditions:
\begin{enumerate}[nolistsep, label=(\arabic*)]
	\item The one with fewer passings is smaller.
	\item The one with fewer crossings is smaller.
	\item The one with fewer components is smaller.
	\item The one with fewer segments is smaller.
	\item For each pair of segments with the same index check:
	\begin{enumerate}[nolistsep, label=(\arabic{enumi}.\arabic*)]
		\item If one has no boundary and the other does, the one with boundary belongs to the smaller PGP.
		\item If there is a boundary, the one with the smaller tuple $B_s$ (with respect to lexicographic order) belongs to the smaller PGP.
		\item The one with more crossings belongs to the smaller PGP.
		\item for each pair of crossings with the same index check:
		\begin{enumerate}[nolistsep, label=(\arabic{enumi}.\arabic{enumii}.\arabic*)]
			\item The one with the smaller absolute value belongs to the smaller PGP.
			\item The one with positive sign belongs to the smaller PGP.
			\item The one with positive orientation belongs to the smaller PGP.
		\end{enumerate}
	\end{enumerate}
\end{enumerate}
It is immediate that if none of these conditions are met, the two PGPs are in fact identical. Note also that up to symmetry, every link has a PGP where the first segment starts at the bottom left, and first crossing in that has index 1 and is a positive crossing. We can remove all diagrams that don't satisfy this condition, as we can transform them with a symmetry move to a diagram that satisfies these conditions and will thus be smaller.

We only keep such PGPs which are locally minimal among all PGPs representing the same link, in the following sense: For a given diagram, it is straightforward to find a minimal PGP with that diagram. Any diagram that admits a single simplifying $\Omega_1$, $\Omega_2$ or $\Theta_2$ move cannot have a minimal PGP for that link and thus can be dropped. We also check whether there are $\Omega_3$ or $\Theta_1$ moves that lead to a smaller PGP, and if so we remove those links as well. Finally, we can explore all diagrams that lead to the same annulus diagram: We can generate all periodic diagrams that can be reached from a given diagram via a sequence of $\Theta_1$ moves, and if any of those diagrams admits a simplifying move, it is not locally minimal. Here we also allow a slightly more general move, the \emph{tangle $\Theta_1$} move, also called $\overline{\Theta}_1$ (see \cref{genMoves}).

We perform this algorithm for links with up to 6 crossings, i.e. $n=6$. The initial number of possible basic diagrams with up to 6 crossings and up to 12 passings is well over 20 million. After performing the minimization as described here and removing non-minimal diagram, we are left with just 2474 diagrams, which we need to check for equivalence.
\subsection{Checking for ambient isotopy equivalence and symmetry}
At this point, for all the remaining links in our list, we compute the relevant invariants: the Alexander and Kauffman polynomials, as well as the hyperbolic volume. To calculate the volume, we use SnapPy \cite{snappy}.\footnote{Of course not all links are hyperbolic; we use the convention that non-hyperbolic links have hyperbolic volume equal to zero.} If we have multi-component links, we can additionally compute these invariants for any sublink.

\begin{figure}
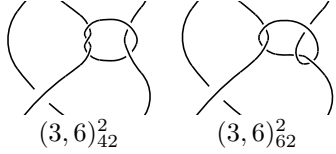

	\begin{center}
		\begin{tikzpicture}[use Hobby shortcut, scale=0.5]
			\input{diagrams/links-1228.tex}
			\node at (2.201, -0.5) {$(3,6)^2_{42}$};
		\end{tikzpicture}
		\begin{tikzpicture}[use Hobby shortcut, scale=0.5]
			\input{diagrams/links-1272.tex}
			\node at (2.305, -0.5) {$(3,6)^2_{62}$};
		\end{tikzpicture}
		\caption{Example of hyperbolic links on which all our invariants coincide, but that can be distinguished by their sublinks: The components of $(3,6)^2_{42}$ are $(2,1)_1^1 \# 3_1^1$ and the longitude knot $(1,0)_1^1$, whereas for $(3,6)^2_{62}$ the components are $(2,1)_1^1$ and $(1,0)_1^1 \# 3_1^1$.}
		\label{linksSameInv}
	\end{center}
\end{figure}

With these, we can partition our set of diagrams into such where all invariants coincide (possibly up to symmetry and orientation change). This gives us 1392 distinct classes with between 1 and 11 elements each. This set contains 51 pairs, 7 triples and three 4-tuples of equivalence classes that have the same Kauffman polynomial, of those there are 33 pairs, 3 triples and one 4-tuple that also have the same Alexander polynomial and hyperbolic volume. Almost all of these cases are composite links with multiple components, where the only difference is along which component the connect sum is taken. The other 7 pairs are prime links; however all of them are multi-component and can be distinguished by considering the sublinks, see \cref{linksSameInv}.

Within these groups, we check whether we can find a sequence of moves that takes one into the other. Due to computational limitations, here we only consider moves that do not increase the number of crossings or the number of passings of the diagram.

\begin{figure}
\begin{center}
	\input{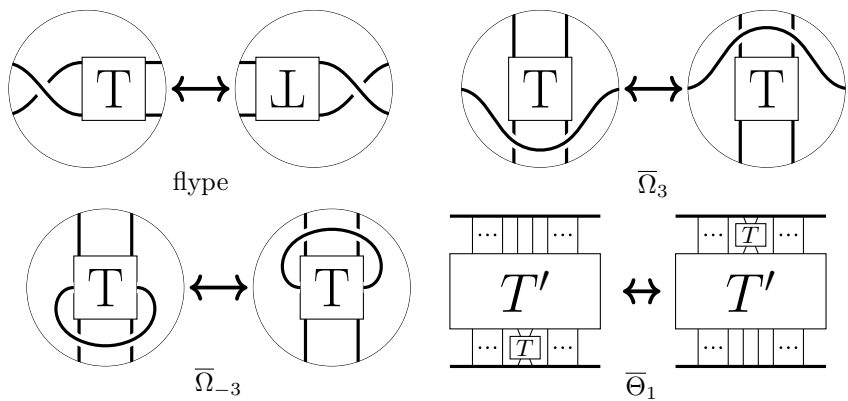}
\end{center}
\caption{More general moves involving tangles.}
\label{genMoves}
\end{figure}

Since this makes our set of moves rather limited, we introduce three new moves: the \emph{flype}, the \emph{$\overline{\Omega}_3$ move} and the \emph{general $\overline{\Omega}_{-3}$ move}, see \cref{genMoves}.
Here, $T$ stands for any tangle. In the flype, we take the mirror of the tangle and flip all crossings, indicated by the upside-down $T$. For the $\overline{\Omega}_3$ and $\overline{\Omega}_{-3}$ move, we will also allow the case where the moved strand has 2 undercrossings rather than 2 overcrossings. It is worth pointing out that in the classical case, the moves $\overline{\Omega}_{\pm 3}$ are equivalent, by simply reversing the roles of the outer and inner tangle (we may view a classical diagram as a diagram in $S^2$, in which case the circles in \cref{genMoves} also bound disks on the outside, which forms a 2-tangle diagram).

These moves do not increase the complexity of the diagram, but if we were to express them as a sequence of $\Omega_i$-moves we might have to increase the number of crossings in between. By defining these as their own moves, we can detect more equivalent links without massively increasing computation times.

The search for a path of moves from one diagram to the other uses a depth-first search algorithm, putting a limit on the number of moves we allow (10 moves are sufficient to find any path like this that can be found). We get 13 cases for which the algorithm fails to find a path; this is because here, we require moves that increase the complexity of the diagram. In these cases, we can find a path of moves by hand. This allows us to conclude that any two links in each of the 1392 classes must be symmetric to one another.

We have now found 1082 links that are symmetric to links with a smaller PGP, and all remaining 1392 links can be distinguished by one of our invariants. Of those, 26 are affine and 163 are non-affine composite links. This means we are left with 526 non-affine prime knots and 658 non-affine prime links with more than one component (a verification that these really are prime can be found in \cref{minprime}). This matches and thus confirms previous findings in \cite{gab12} and \cite{sch16}.
\subsection{Checking for Dehn twists}
We have completed our classification up to a-equivalence (and indeed, up to symmetry). Next, we are concerned with determining the h-equivalence classes. For the remaining symmetry classes of links, we want to find out which of them are the result of Dehn twisting simpler links. Note that any homeomorphism keeps the hyperbolic volume intact, so any links with different hyperbolic volume cannot be h-equivalent.

We can thus sort our links by their volume. This is a first distinction, but it still leads to fairly large classes: The set of non-hyperbolic links alone has 301 elements, but even the set of all hyperbolic links with the smallest volume still contains 201 elements - for these, the volume is always approximately 3.66386..., which is known to be the smallest (nonzero) volume that a complement of any 2-component classical link can have (see \cite{ag10}). It is clear that we still need to do better.

For any two links with the same volume, we can first check whether their wrapping numbers and maximal winding numbers coincide - a necessary condition for homeomorphic equivalence - and then see whether the invariants $\mathcal{A}_L$ coincide. In cases where the maximal winding number $w \neq 0$, and $\deg_x(\overline{\Delta}_L(x,t)) > 0$, we can determine the degree of a Dehn twist $d^k$ by \cref{uniqueTwist}. Otherwise, we can use our results on the Kauffman polynomial to obtain bounds for $k$. There are cases where neither of these methods applies; there we just test twists with $|k| \leq 3$; it turns out that with this, in each of these cases we find a Dehn twist take takes one to the other.

Once we have a possible twist, we generate the twisted diagram, and subsequently we may calculate the Kauffman polynomial of the diagram to see whether it matches up with the Kauffman polynomial of the other link. If not, we already know that this is not a valid option. Otherwise, we follow a depth-first search to find a Reidemeister path from one link to the other.

Finally, it might also be the case that a link $L$ is a Dehn twist not of another link class in our list, but instead it is h-split, i.e. the Dehn twist of a split link consisting of multiple links form our classes. We use the next proposition to rule out this possibility for most of our classes. First, we prove a folklore statement about the classical multivariable Alexander polynomial $\Delta$.
\begin{lemma}
	\label{connAlex}
	Let $L_i \subseteq S^3$ be a $k_i$-component oriented classical link with $k_i > 1$, for $i \in \{1,2\}$, and let $L_1 \# L_2$ be a connected sum along the first components. Then
	\begin{equation}
		\Delta_{L_1 \# L_2}(r, t_2, \ldots, t_{k_1}, s_2, \ldots, s_{k_2}) = (r-1)\Delta_{L_1}(r, t_2, \ldots, t_{k_1}) \Delta_{L_2}(r, s_2, \ldots, s_{k_2}).
	\end{equation}
	Here the first component of $L_i$ has its meridian evaluating to $r$.
\end{lemma}
\begin{proof}
	We determine the fundamental group of the complement of the connected sum $L_1 \# L_2$ via Wirtinger presentation. Consider a link diagram of $L_1 \# L_2$ that exhibits the connected sum, i.e. there is a circle $C$ that intersects the diagram exactly twice (and the diagram is connected).
	
	Now consider the generating loops of the Wirtinger presentation around each arc of the diagram. We can single out the two loops $s_0^1, s_0^2$ around the arcs $a_1, a_2$ that intersect $C$. Note that these loops are actually homotopic within the link complement, and thus represent the same generator within $\pi_1 (S^3 \setminus L_1 \# L_2)$. Comparing the resulting presentation to the presentations obtained for $L_1$ and $L_2$ by considering their diagrams obtained by undoing the connected sum, it is clear that, given the resulting presentations $\left\langle s_0^1, \ldots s_{k_1}^1 \, \middle| \, r_0^1, \ldots, r_{l_1}^1 \right\rangle$ and $\left\langle s_0^2, \ldots s_{k_2}^2 \, \middle| \, r_0^2, \ldots, r_{l_2}^2 \right\rangle$, we can write the presentation for $\pi_1 (S^3 \setminus L_1 \# L_2)$ as
	\begin{equation}
		\left\langle t_0, \ldots t_{k_1}, s_0, \ldots s_{k_2} \, \middle| \, r_0^1, \ldots, r_{l_1}^1, r_0^2, \ldots, r_{l_2}^2,  s_0(t_0)^{-1} \right\rangle
	\end{equation}
	We may substitute $t_0$ and $s_0$ by $r$ to simplify the presentation. Then we do Fox calculus (for details on this method of computing the Alexander polynomial, see \cite{fox}) to obtain the associated matrix of free derivatives. The relations split into 2 sets that only have the generator $r_0$ in common, and thus the matrix obtained is almost block diagonal, with the exception of the column belonging to $r_0$. To get the Alexander polynomial, we can strike that column, and the resulting matrix $A$ is block diagonal with blocks $A_1$, $A_2$ that correspond to the matrices we get from doing the same to the individual presentations for $L_1$, $L_2$. Thus $\det (A) = \det(A_1) \det(A_2)$, and since $\Delta_{L_i} \doteq \frac{\det (A_i)}{r-1}$, the claim follows.
\end{proof}
\begin{prop}
	\label{splitAlex}
	Assume that an oriented link $L \subseteq V^3$ is homeomorphically equivalent to a split link $L_1 | L_2$. Then we have
	\begin{equation}
	\Delta_L(t,x) \doteq (xt^{wk} - 1) \Delta_{L_1}(t, xt^{wk}) \Delta_{L_2} (t, xt^{wk}),
	\end{equation}
	where $w$ is the winding number of $L_1|L_2$, for some $k \in \Z$, and $L$ does not have a hyperbolic complement.
\end{prop}
\begin{figure}
	\begin{center}
		\input{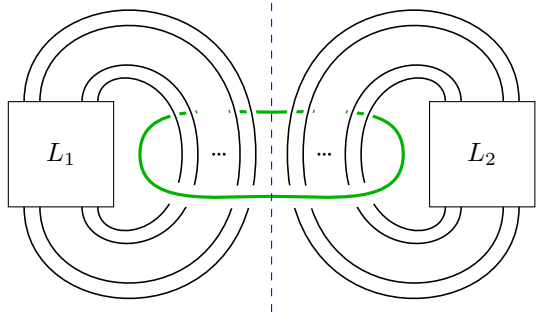}
		\caption{The augmentation of a split link is composite.}
		\label{splitAug}
	\end{center}
\end{figure}
\begin{proof}
	We can see that the augmentation link of $L_1|L_2$ is the connect sum $\au(L_1) \# \au(L_2)$ along the augmentation component, see \cref{splitAug}. As a composite link, $\au(L_1|L_2)$ is non-hyperbolic, and thus neither is $L$.
	
	By \cref{connAlex} we know that the Alexander polynomial is
	\begin{equation}
		\begin{split}
		\Delta_{L_1|L_2}(t,x) \doteq \Delta_{\au(L_1|L_2)}(t,x) \doteq (x-1) \Delta_{\au(L_1)}(t,x) \Delta_{\au(L_2)}(t,x) \\ \doteq (x-1)\Delta_{L_1}(t,x)\Delta_{L_2}(t,x).
		\end{split}
	\end{equation}
	But from \cref{bicAlex} we then know that $\Delta_L(t,x) \doteq \Delta_{L_1|L_2}(t,x^{wk})$ for some $k \in \Z$, and thus the claim follows.
\end{proof}
This means that for non-hyperbolic links in the solid torus, we can check whether the polynomials in $\mathcal{A}_L$ have a factorization into Alexander polynomials and $(xt^k-1)$. If a factorization exists, we find all combinations of non-split links whose combinations might form something h-equivalent to our original link. Necessary conditions are that the total winding numbers sum up to the total winding number of the link in question, and analogously for the wrapping number. Furthermore, the Alexander polynomials of the split summands must appear as factors in the Alexander polynomial of their combination. For all valid combinations, we check whether they are Dehn twists as before.

\begin{figure}
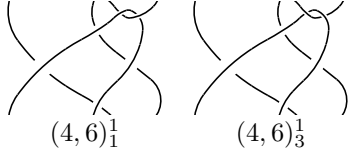

	\begin{center}
		\begin{tikzpicture}[use Hobby shortcut, scale=0.5]
			\input{tikz/links-1535.tex}
			\node at (2.357, -0.5) {$(4,6)^1_{1}$};
		\end{tikzpicture}
		\begin{tikzpicture}[use Hobby shortcut, scale=0.5]
			\input{tikz/links-1539.tex}
			\node at (2.358, -0.5) {$(4,6)^1_{3}$};
		\end{tikzpicture}
		\caption{Solid torus knots for which our algorithm fails to determine h-equivalence.}
		\label{specialLinks}
	\end{center}
\end{figure}

These methods allow us to decide the question of h-equivalence for almost all links in our list. Only a single pair of link classes remains for which the algorithm does not yield a definitive result, see \cref{specialLinks}. We will provide a proof here that these two knots are not h-equivalent.

We will need the following notion, which is due to Schubert \cite{sc49}.
\begin{definition}
	Let $K$ be a knot in a 3-manifold $M$. $K$ is called a \emph{satellite knot} if the complement $M \setminus K$ contains an incompressible, non-boundary-parallel torus.
\end{definition}
We will skip the definitions of incompressible and boundary parallel surfaces for brevity; the interested reader may find them in e.g. \cite[Chapter 3]{sc14}. We have the following observation:
\begin{lemma}
	$K \subset V^3$ is a satellite knot if and only if its Dehn twist $d(K)$ is a satellite knot.
\end{lemma}
\begin{proof}
	The complements of $K$ and $d(K)$ are homeomorphic via the restriction of $d$, and the property of containing an incompressible, non-boundary-parallel torus is preserved under homeomorphisms.
\end{proof}
\begin{prop}
	The knots $K_1 = (4,6)^1_{1}$ and $K_2 = (4,6)^1_{3}$ as seen in \cref{specialLinks} are not h-equivalent.
\end{prop}
\begin{proof}
	First, note the Kauffman polynomials of the two knots:
	\begin{align*}
		\nabla_{K_1}(A,x) = & \left(A^{10} - A^6\right)x^4 + \left( -3A^{10}+3A^6-A^2+A^{-2}\right)x^2 - A^6 - A^{-6} \\
		\nabla_{K_2}(A,x) = & \left(-A^2+A^{-2}\right)x^4 + \left(3A^2-3A^{-2}+A^{-6}-A^{-10}\right)x^2 \\ &- A^2 - A^{-6} - A^{-10} + A^{-14}
	\end{align*}
	This immediately gives that $K_1, K_2$ are not symmetric. Furthermore we can compute
	\begin{multline*}
		\nabla_{d^{-1}(K_1)}(A,x) = x^4(-A^{-22} + A^{26}) + x^2(A^{-14} = A^{-18} + 3A^{-22} - 3 A^{-26}) \\- A^{-2} - 2A^{-10} + A^{-18} - A^{-22} + A^{26},
	\end{multline*}
	so $d^{-1}(K_1)$ and $K_2$ are not symmetric. We will next make an argument that $d^k(K_1)$ and $K_2$ are not symmetric for any $k \notin \{ 0,-1\}$. To this end, note that $K_1, K_2$ are both satellite knots, with the incompressible torus given as the boundary of a regular neighborhood $N$ of the knot $K_o = (2,1)^1_1$.
	
	In the solid torus $N$ with its obvious marking, $K_1$ forms the knot $K_i = (2,2)^1_1$. Consider now a $k$-fold Dehn twist of $K_1$; this will map $N$ to $d^k(N)$, a regular neighborhood of $d^k(K_o)$; and within $d^k(N)$, $d^k(K_1)$ will form the link $d^{2k}(K_i)$. The Dehn twist in within $d^k(N)$ is doubled since $K_o$ has winding number 2.
	
	Now for any knot $K$, we may consider the augmentation link $\au(K) \subset S^3$, and then obtain a knot in $S^3$ by removing the augmentation component; let us call the resulting knot $\au_{\red}(K)$, its embedding into the round solid torus in $S^3$. This is clearly also an a-equivalence invariant. It is not hard to see that $\au_{\red}(K_2) = \au_{\red}(K_1)$ is the unknot. If we can prove that $\au_{\red}(d^k(K_1))$ is not the unknot, then this would finish the proof.
	
	We note that $\au_{\red}(d^k(K_o))$ is the $(2,2k+1)$-torus knot and thus - since $k \neq 0,-1$ -- not the unknot (this is a standard result and may be found e.g. in \cite[Theorem 2.2.2]{ka12}). Furthermore, $d^{2k}(K_i)$ is not affine in $d^k(N)$ as we have already determined via the methods in this chapter that $K_i$ is not h-equivalent to an affine knot. In particular, this implies that $\au_{\red}(d^k(K_1))$ is a satellite knot, and finally a theorem of Rolfsen \cite[Corollary D.10]{ro90} says that in particular $\au_{\red}(d^k(K_1))$ is nontrivial.
\end{proof}
This settles the last remaining case. In total we find 105 cases of symmetry classes of links where either the h-equivalence class is generated by a symmetry class with smaller minimal PGP, or where the symmetry class is h-equivalent to a split link consisting of two or more other links in our list (i.e. with 6 or fewer crossings). That means we arrive at 1287 different h-equivalence classes for links in the solid torus with up to 6 crossings.

\subsection{Generating diagrams and link tables}
In order to get nice expressive diagrams for these links, we use a circle packing algorithm \cite{cs03}. Given a periodic Gauss paragraph, we consider the underlying annulus graph $\Gamma$. We may identify both boundary components of the annulus to obtain a graph embedded in the torus. To this, we add the image of the annulus boundary as a cyclic graph $C$ with two vertices. Then we apply barycentric subdivision to the graph $\Gamma \cup C$, and thus obtain a triangulation of the torus $T^2$. Now, there is a unique (up to scaling) circle packing associated with this triangulation (i.e. each there is a circle for each triangle, and two circles touch if and only if their corresponding triangles share an edge), which determines a flat metric on $T^2$. For more background on this method, a reader may want to check \cite{BS90}.

We can use this circle packing, which has a circle for every crossing, edge and face of the periodic link diagram, and draw the link through the circles associated to crossings and edges. The circle packing guarantees that we won't get any unintended crossings of the drawn lines.

We also choose a naming convention for our link classes: Any prime link will have a a name of the form $(\wrap, k)^c_i$. Here, $\wrap$ is the wrapping number, $k$ is the crossing number, $c$ is the number of components of the link, and $i$ is a running index.
\begin{figure}
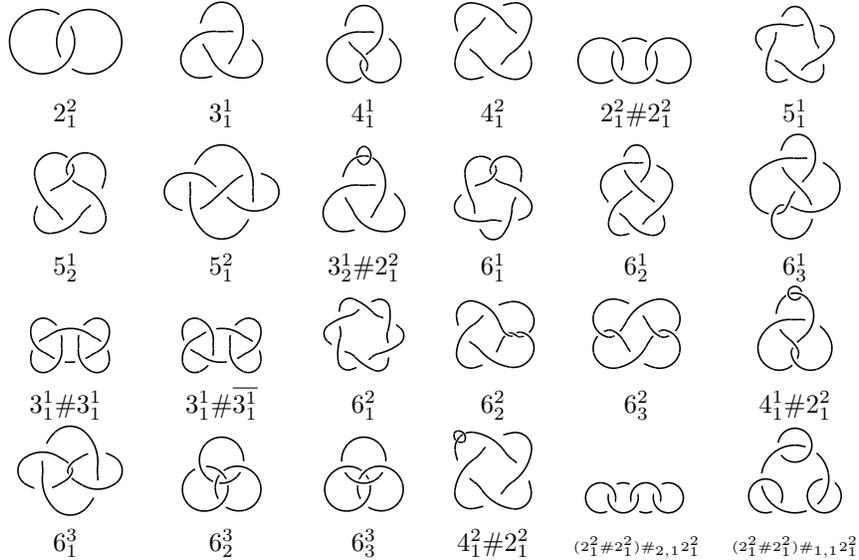

	\begin{tabular}{c c c c c c}
		\begin{tikzpicture}[use Hobby shortcut, scale=0.3]
			\input{diagrams/links-2.tex}
		\end{tikzpicture} &
		\begin{tikzpicture}[use Hobby shortcut, scale=0.3]
			\input{diagrams/links-3.tex}
		\end{tikzpicture} &
		\begin{tikzpicture}[use Hobby shortcut, scale=0.3]
			\input{diagrams/links-4.tex}
		\end{tikzpicture} &
		\begin{tikzpicture}[use Hobby shortcut, scale=0.3]
			\input{diagrams/links-5.tex}
		\end{tikzpicture} &
		\begin{tikzpicture}[use Hobby shortcut, scale=0.3]
			\input{diagrams/links-6.tex}
		\end{tikzpicture} &
		\begin{tikzpicture}[use Hobby shortcut, scale=0.3]
			\input{diagrams/links-8.tex}
		\end{tikzpicture} \\
		$2_1^2$ & $3_1^1$ & $4_1^1$ & $4_1^2$ & $2_1^2 \# 2_1^2$ & $5_1^1$ \\ \rule{0pt}{0.5cm}
		\begin{tikzpicture}[use Hobby shortcut, scale=0.3]
			\input{diagrams/links-7.tex}
		\end{tikzpicture} &
		\begin{tikzpicture}[use Hobby shortcut, scale=0.3]
			\input{diagrams/links-11.tex}
		\end{tikzpicture} &
		\begin{tikzpicture}[use Hobby shortcut, scale=0.3]
			\input{diagrams/links-12.tex}
		\end{tikzpicture} &
		\begin{tikzpicture}[use Hobby shortcut, scale=0.3]
			\input{diagrams/links-19.tex}
		\end{tikzpicture} &
		\begin{tikzpicture}[use Hobby shortcut, scale=0.3]
			\input{diagrams/links-20.tex}
		\end{tikzpicture} &
		\begin{tikzpicture}[use Hobby shortcut, scale=0.3]
			\input{diagrams/links-16.tex}
		\end{tikzpicture} \\
		 $5_2^1$ & $5_1^2$ & $3^1_2 \# 2^2_1$ & $6^1_{1}$ & $6^1_{2}$ & $6^1_{3}$ \\ \rule{0pt}{0.5cm}
		 \begin{tikzpicture}[use Hobby shortcut, scale=0.2]
		 	\input{diagrams/links-13.tex}
		 \end{tikzpicture} &
		 \begin{tikzpicture}[use Hobby shortcut, scale=0.2]
		 	\input{diagrams/links-15.tex}
		 \end{tikzpicture} &
		 \begin{tikzpicture}[use Hobby shortcut, scale=0.3]
		 	\input{diagrams/links-36.tex}
		 \end{tikzpicture} &
		 \begin{tikzpicture}[use Hobby shortcut, scale=0.3]
		 	\input{diagrams/links-32.tex}
		 \end{tikzpicture} &
		\begin{tikzpicture}[use Hobby shortcut, scale=0.3]
			\input{diagrams/links-34.tex}
		\end{tikzpicture} &
		\begin{tikzpicture}[use Hobby shortcut, scale=0.3]
			\input{diagrams/links-41.tex}
		\end{tikzpicture} \\
		$3^1_1 \# 3^1_1$ & $3^1_1 \# \overline{3^1_1}$ & $6^2_1$ &$6^2_{2}$ & $6^2_{3}$ & $4^1_1 \# 2^2_1$ \\ \rule{0pt}{0.5cm}
		\begin{tikzpicture}[use Hobby shortcut, scale=0.3]
			\input{diagrams/links-43.tex}
		\end{tikzpicture} &
		\begin{tikzpicture}[use Hobby shortcut, scale=0.3]
			\input{diagrams/links-51.tex}
		\end{tikzpicture} &
		\begin{tikzpicture}[use Hobby shortcut, scale=0.3]
			\input{diagrams/links-42.tex}
		\end{tikzpicture} &
		\begin{tikzpicture}[use Hobby shortcut, scale=0.3]
			\input{diagrams/links-45.tex}
		\end{tikzpicture} &
		\begin{tikzpicture}[use Hobby shortcut, scale=0.2]
			\input{diagrams/links-55.tex}
		\end{tikzpicture} &
		\begin{tikzpicture}[use Hobby shortcut, scale=0.3]
			\input{diagrams/links-57.tex}
		\end{tikzpicture} \\
		$6^3_{1}$ & $6^3_{2}$ & $6^3_{3}$ & $4^2_1 \# 2^2_1$ & {\tiny $(2^2_1 \# 2^2_1) \#_{2,1} 2^2_1$} & {\tiny $(2^2_1 \# 2^2_1) \#_{1,1} 2^2_1$}
	\end{tabular}
	\caption{The classical non-split links up to 6 crossings, up to symmetry (excluding the empty link and the unknot)}
	\label{affineLinks}
\end{figure}
Composite links are named as $N_1 \# N_2$, where $N_1$ and $N_2$ are the names of its summands with $N_2$ the name of the affine summand (which we identify with the corresponding link in $S^3$). If the affine component is the mirror image of one of the links from \cref{affineLinks}, we write an overline to indicate this, e.g. $\overline{3_1^1}$ for the mirror of the trefoil. We use brackets in cases of composite links with connect summands. Note that this naming convention need not produce unique names in cases where $N_1$ or $N_2$ have more than one component. In these instances, we improve our notation by subscripts $K \#_{i,j} L$ where $i$ denotes the segment along which we connect sum $L$ (in the order they appear in its minimal PGP), and $j$ the segment along which we connect sum $K$.

Now we can collect all our information to give a link table of h-equivalence classes for solid torus links. Only one question remains to be addressed.
\subsection{Minimality, primeness and splitness}
\label{minprime}
Our link classes are labeled by wrapping number and crossing number, which are obtained from their diagrams. In order for this to make sense, we need to be sure that the diagrams that realize the crossing number rsp. wrapping number of a given link class are indeed among the diagrams we generated. For instance, it may be the case that we generated a 5-crossing diagram with wrapping number 3, but there is a diagram of the same link (up to h-equivalence) with 17 crossings and wrapping number 2.

Showing that we found the minimal crossing number is simple enough, since we created all possible diagrams with lower crossing numbers and saw that any given link in our list is not h-equivalent to any of the others. To show that our wrapping numbers are minimal, we use \cref{propWrappingBound}. In fact, all links in our table satisfy \cref{conjWrapping}, which proves that the wrapping numbers we found must be minimal.

We also distinguish between prime and non-prime links. For this to make sense, we need to show that the link classes we claim to be prime indeed cannot have composite representatives. For this, first note that primeness is preserved under homeomorphism (as the ball in which we form the connected sum is mapped to a ball), so either all representatives are prime or none. No composite link can be hyperbolic \cite{thurnotes}, so we can restrict the problem to only the non-hyperbolic cases. This already reduces the problem significantly: There are only 230 non-hyperbolic links in our list. Of those, 26 are affine (and primeness is solved for affine links), and of the remaining 204 we have 160 in which the diagram already is composite. Only in 44 cases do we need to apply additional methods.

We rely on the fact that in the case of classical knots and links, we already have tables of all prime links up to 14 crossings (and indeed much more for knots) available in SnapPy; compare \cite{ho05} and \cite{snappy}. The next proposition allows us to reduce the problem to the classical case via the augmentation.
\begin{prop}
	Assume that $L = L_r \# A$ is a composite link in the solid torus. Then $\au(L)$ is composite as well.
\end{prop}
\begin{proof}
	If $L$ is composite, then there is a 2-sphere $S \subset V^3$ that intersects $L$ twice and splits it into two nontrivial components. The image of this sphere under the embedding $g: V^3 \to S^3$ splits $\au(L)$ into two components, one of them $A$ (which is not the unknot), and the other containing the augmentation component, so neither can be trivial.
\end{proof}
Unfortunately, this does not easily work in the other direction; a knotting for which $\au(K)$ is composite may still be prime (for instance this is the case for any split monoperiodic knotting). We do however gain the result that if $\au(K)$ is prime, so is $K$. We find that for all the 44 links that we had left, the augmentation is prime, thus proving that the links themselves are prime.

\begin{figure}
	\begin{center}
	\input{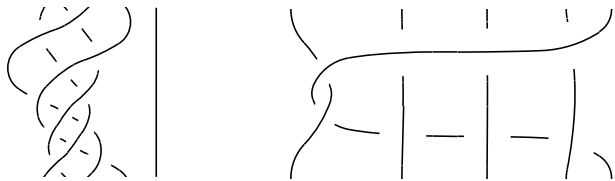}
	\caption{Example of a split link which is h-equivalent to a link with lower crossing number.}
	\label{splitMin}
	\end{center}
\end{figure}
A final claim that we would like to make is that all links we present are indeed non-split up to h-equivalence. Unfortunately it is not the case that the minimal diagram up to h-equivalence of a split link is the split diagram, as evidenced by \cref{splitMin}. Since the augmentation of a split link is composite, we can use techniques analogous to the case of primeness: Hyperbolic links cannot be h-split, and a composite link is h-split if and only if its non-affine summand is. Furthermore, clearly knots can never be h-split.

We can thus limit our investigation to prime, non-hyperbolic links with at least 2 components. Our list contains 50 such cases. Additionally, any split link must have a factor of the form $(xt^k - 1)$ in its Alexander polynomial, where $k$ is a multiple of its total winding number, by \cref{splitAlex}. Only the two cases in the following proposition satisfy this condition.
\begin{prop}
	The following solid torus links are not h-equivalent to split links:
	
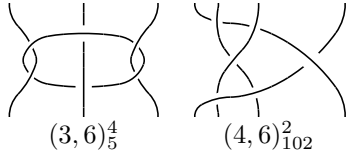
\begin{center}
		\begin{tikzpicture}[use Hobby shortcut, scale=0.5]
			\input{diagrams/links-1428.tex}
			\node at (2.299, -0.5) {$(3,6)^4_{5}$};
		\end{tikzpicture}
		\hspace{2cm}
		\begin{tikzpicture}[use Hobby shortcut, scale=0.5]
			\input{diagrams/links-1884.tex}
			\node at (2.368, -0.5) {$(4,6)^2_{102}$};
		\end{tikzpicture}
	\end{center}
\end{prop}
\begin{proof}
	First, let us look at an arbitrary split link $L = L_1|L_2$, and let $g: V^3 \to S^3$ be the marking embedding. Then for any two components $K_1 \subseteq L_1$ and $K_2 \subseteq L_2$, it is clear that $\lk\big(g(K_1), g(K_2)\big) = 0$ (as there is a diagram of $\au(L)$ in which they share no crossing). We can use this fact to prove that the links in question are not split, and indeed not -hequivalent to split links.
	
	In the case of $(3,6)^4_{5}$, we note that all the non-affine components have nonzero linking number with the affine component. The linking number with an affine component remains the same (up to sign) under homeomorphisms, so there is no way to partition the components of $(3,6)^4_{5}$ into two nonemtpy subsets where each two elements of different subsets have linking number zero. By the preceding observation, this implies that $(3,6)^4_{5}$ cannot be split.
	
	We can compute that the components of $(4,6)^2_{102}$ have linking number 0 with one another (in the augmentation); a Dehn twist changes that by $\pm 4$ (the product of the winding numbers), so for each nontrivial $k$-fold Dehn twist we get that the link components have nonzero linking number, and the link cannot be split. The only representative of the h-equivalence class that could still be split is $(4,6)^2_{102}$ itself.
	
	Assume that $(4,6)^2_{102}$ was split; then it must consist of its two components, both of the form $(2,1)^1_1$. But if we compare the Alexander polynomial of the split sum $(2,1)^1_1|(2,1)^1_1$ with that of $(4,6)^2_{102}$, we see that they differ; thus the link cannot be split.
\end{proof}

\section{Results}
We find that of the 1366 non-affine non-split knots and links in the solid torus with crossing number up to 6 up to symmetry, 163 are composite. 
\begin{figure}
	\input{tikz/dehn-removed.tex}
	\caption{Prime diagrams that were found h-equivalent to simpler ones and do not immediately admit a $\Delta$ move. The twisted diagram are in the top row, their simpler versions in the bottom row (with the name we gave that h-equivalence class). For ease of notation, we write $\ell = (1,0)_1^1$ for the longitude knot.}
\end{figure}
Of the prime knots, we find that there are 37 that can be reduced to simpler knots via h-equivalence. For the prime links, there are 39 cases that are h-equivalent either to simpler links or to split links. Of the composite knots and links, 29 can be reduced in this way. Thus, in total our link table up to h-equivalence contains 1261 non-affine non-split knots and links in the solid torus.

The classification can serve as a good reference to come up with new theorems. Let us give an example. From our link table, it seems that the links with the highest wrapping number (and fixed crossing number) are always knots. Indeed we can show:
\begin{prop}
	\label{windProp}
	Let $L$ be a non-split link in $V^3$ with $n$ components, maximal winding number $w$ and crossing number $c$. Then
	\begin{equation}
		\label{windBound}
		w \leq c - n + 2.
	\end{equation}
\end{prop}
\begin{proof}
	Note that we showed in \cref{passBound} that the wrapping number $\wrap$ of $L$ is at most $c+1$. If $L$ is a knot (i.e. $n=1$), we can conclude from this that
	\begin{equation}
		w \leq \wrap \leq c+1 = c - n + 2.
	\end{equation}
	Now suppose we have shown the statement for any $k < n$. We may write $L$ as the union of two non-empty sublinks $L = L_1 \cup L_2$. Let $n_i, w_i$ and $c_i$ be the number of components, winding number and crossing number of $L_i$, respectively.
	
	If we have a diagram $D$ of $L$ with $c$ crossings, we may extract from this diagrams $D_i$ of the $L_i$ by removing the other component. For $L$ to be non-split, $L_1$ and $L_2$ need to be linked, i.e. there need to be at least 2 crossings in $D$ where $L_1$ and $L_2$ meet. These crossings will not be present in the diagrams $D_i$, hence $D_1$ and $D_2$ together will have at most $c - 2$ crossings. From this follows that $c_1 + c_2 \leq c - 2$.
	
	It is also clear from the definition of the total winding number that $w = w_1 + w_2$. Finally, since $n_i < n$, we have that the statement holds for the $L_i$, i.e.
	\begin{equation}
		w_i \leq c_i - n_i + 2.
	\end{equation}
	But then it follows that
	\begin{equation}
		w = w_1 + w_2 \leq c_1 + c_2 - n_1 + n_2 + 4 \leq c-n + 2.
	\end{equation}
\end{proof}
Indeed our results suggest that \cref{windBound} holds not only for the winding number, but indeed for the wrapping number.

We can also note that in our tabulation, every link has a diagram which realizes both the crossing number and the wrapping number at the same time. It seems reasonable to expect that this is generally possible.
\begin{conjecture}
	\label{wrapConject}
	Let $L$ be a solid torus link with crossing number $c$ and wrapping number $\wrap$. Then there is a diagram $D$ of $L$ which has $c$ crossings and $\wrap$ passings.
\end{conjecture}
If this is true, we could drop the technical assumption from our next proposition, thereby making a more general statement about the relation between wrapping and winding number.
\begin{prop}
	Let $L$ be a non-split solid torus link with crossing number $c$ and wrapping number $\wrap = c+1$. Assume that there is a diagram $D$ for $L$ with $c$ crossings and $\wrap$ passings. Then $L$ is a knot, and $c+1$ is also the winding number of $L$.
\end{prop}
\begin{proof}
	If we can show that the winding number $w$ of $L$ is equal to the wrapping number (i.e. equal to $c+1$), then the other part of the statement follows from \cref{windProp}. So let us assume towards a contradiction that $w < \wrap = c+1$.
	
	We will show towards a contradiction that $D$ is equivalent to a diagram with fewer passings. For this, form the graph $\Gamma_e$ with a vertex or each segment of $D$ and and edge between two vertices for each crossing that the corresponding segments share. First let us assume that $\Gamma_e$ is disconnected. This is equivalent to $D$ being disconnected and splits into two parts $D = D_1 \cup D_2$. By assumption, $L$ is non-split, so $D$ is not split, thus $D_1,D_2$ cannot be valid diagrams, which means they have a different number of intersections with the boundary to either side. 
	
	But then, there is a sloped line separating $D_1$ from $D_2$, and we can find a cutting curve for $D$ that follows this sloped line for a time such that the cutting curve intersects $D_i$ on the side where it has fewer intersections with the boundary. Transforming the diagram with this cutting curve gives a diagram with lower wrapping number, a contradiction. Thus $D$ (and also $\Gamma_e$) are connected.
	
	Note the $\Gamma_e$ has $c$ edges and $c+1$ vertices, so it forms a tree. Hence no segment in $D$ has any self-crossing (which would give a loop in $\Gamma_e$), and no two segments of $D$ can cross twice - i.e. $D$ is a basic diagram.
	
	If each segment of $D$ started and ended at opposite edges of the square, orienting them going from top to bottom would yield a valid orientation for the link, and with this we could compute that $w = \wrap$, which contradicts our assumption. Thus there exists a segment $s_0$ that starts and ends at the same edge.
	
	We consider the region $R$ in $D$ bounded by $s_0$ and the edge it has its endpoints on. We may assume that $R$ contains no other segment completely; otherwise we take $s_0$ to be that smaller segment. Now if $R$ intersects no other segments, that implies that $s_0$ has no crossings and so we can perform a $\Theta_2$ move lowering the number of passings; a contradiction.
	
	Let $t_1, \ldots, t_k$ be the segments intersecting $R$ (and thus crossing $s_0$). Note that $t_i$ and $t_j$ do not cross for any $i,j$; if they did, then $t_i, t_j, s_0$ would from a cycle in $\Gamma_e$. Thus we can use a $\Theta_1$ move on the outermost crossing of $s_0$. The resulting diagram $\widetilde{D}_1$ still has $c$ crossings and $c+1$ passings, and now the corresponding region $\widetilde{R}_1$ intersects one segment less. Continuing inductively, we obtain a diagram $\widetilde{D}_k$ such that $\widetilde{R}_k$ intersects no other segments, which is a contradiction.
\end{proof}

\section{Table of solid torus links up to 6 crossings}
As mentioned previously, all knots and links with wrapping number less or equal to 1 can be extracted from the affine link table (\cref{affineLinks}), so we shall start our table at wrapping number 2. We present each h-equivalence class of non-split links in the solid torus with up to 6 crossings in a periodic diagram, ordered by wrapping number first, crossing numer second, and number of components last. Links where all these number coincide are ordered by the order in PGPs established in \cref{secGenOrder}.

The associated invariants (Kauffman and Alexander polynomials as well as hyperbolic volumes) are available in the distribution of the code at \cite{code}.
\begin{center}
\input{strips-in.tex}
\end{center}
\printbibliography
\end{document}

%% file: tikz/skein-start.tex
\begin{tikzpicture}[baseline=2pt, scale=0.15]
	\draw [dotted] (0.5,1) circle (1cm);
	\clip (0.5,1) circle (1cm);
	\draw [white, double=black] (1,0) -- (0,2);
	\draw [white, double=black, line width=1.5pt] (0,0) -- (1,2);
\end{tikzpicture}

%% file: tikz/skein-plus.tex
\begin{tikzpicture}[baseline=2pt, scale=0.15]
	\draw [dotted] (0.5,1) circle (1cm);
	\clip (0.5,1) circle (1cm);
	\draw [white, double=black](0,0) -- (0,2);
	\draw [white, double=black] (1,0) -- (1,2);
\end{tikzpicture}

%% file: tikz/skein-minus.tex
\begin{tikzpicture}[baseline=2pt, scale=0.15]
	\draw [dotted] (0.5,1) circle (1cm);
	\clip (0.5,1) circle (1cm);
	\draw [white, double=black](0.5,0) circle (0.5cm);
	\draw [white, double=black] (0.5,2) circle (0.5cm);
\end{tikzpicture}

%% file: tikz/move-omega-kauff-1.tex
\begin{tikzpicture}[baseline=2pt, scale=0.15, use Hobby shortcut]
	\draw [dotted] (0.5,1) circle (1cm);
	\clip (0.5,1) circle (1cm);
	\draw [white, double=black] ([in angle=90, out angle=90]0,0) .. (1,1);
	\draw [white, double=black, line width=1.5pt] ([in angle=-90, out angle=-90]1,1) .. (0,2);
\end{tikzpicture}

%% file: tikz/move-omega-kauff-2.tex
\begin{tikzpicture}[baseline=2pt, scale=0.15]
	\draw [dotted] (0,1) circle (1cm);
	\clip (0,1) circle (1cm);
	\draw [white, double=black] (0,0) -- (0,2);
\end{tikzpicture}

%% file: tikz/move-omega-kauff-3.tex
\begin{tikzpicture}[baseline=2pt, scale=0.15, use Hobby shortcut, xscale=-1]
	\draw [dotted] (0.5,1) circle (1cm);
	\clip (0.5,1) circle (1cm);
	\draw [white, double=black] ([in angle=90, out angle=90]0,0) .. (1,1);
	\draw [white, double=black, line width=1.5pt] ([in angle=-90, out angle=-90]1,1) .. (0,2);
\end{tikzpicture}

%% file: tikz/skein-invert.tex
\begin{tikzpicture}[baseline=2pt, scale=0.15]
	\draw [dotted] (0.5,1) circle (1cm);
	\clip (0.5,1) circle (1cm);
	\draw [white, double=black] (0,0) -- (1,2);
	\draw [white, double=black, line width=1.5pt] (1,0) -- (0,2);
\end{tikzpicture}

%% file: tikz/kauff-dehn-2-proof-1.tex
\begin{tikzpicture}[use Hobby shortcut, baseline=6pt, scale=0.3]
	\foreach \i in {0,1} {
		\draw ([in angle=-90, out angle=90]1,\i) .. (0,1 + \i);
		\draw [white, double=black, line width=2pt] ([in angle=-90, out angle=90]0,\i) .. (1,1 + \i); }
\end{tikzpicture}

%% file: tikz/kauff-dehn-2-proof-2.tex
\begin{tikzpicture}[use Hobby shortcut, baseline=6pt, scale=0.3]
	\draw ([in angle=-90, out angle=90]1,0) .. (0,2);
	\draw [white, double=black, line width=2pt] ([in angle=-90, out angle=90]0,0) .. (1,2);
\end{tikzpicture}

%% file: tikz/kauff-dehn-2-proof-3.tex
\begin{tikzpicture}[use Hobby shortcut, baseline=6pt, scale=0.3]
	\draw ([in angle=-90, out angle=90]1,1.2) .. (0,2);
	\draw [white, double=black, line width=2pt] ([in angle=-90, out angle=90]0,1.2) .. (1,2);
	\draw ([in angle=90, out angle=90]0,0) .. (1,0);
	\draw ([in angle=-90, out angle=-90]0,1.2) .. (1,1.2);
\end{tikzpicture}

%% file: tikz/kauff-dehn-3-basic-1.tex
\begin{tikzpicture}[scale=0.3, baseline=1pt]
	\clip (0,0) rectangle (0.8,1);
	\draw (0,0) -- (0.8,1);
	\draw (0.6,0) circle (0.2cm);
	\draw (0.2, 1) circle (0.2cm);
\end{tikzpicture}

%% file: tikz/kauff-dehn-3-basic-2.tex
\begin{tikzpicture}[scale=0.3, baseline=1pt]
	\clip (0,0) rectangle (0.8,1);
	\draw (0,1) -- (0.8,0);
	\draw (0.6,1) circle (0.2cm);
	\draw (0.2, 0) circle (0.2cm);
\end{tikzpicture}

%% file: tikz/kauff-dehn-3-proof-1.tex
\begin{tikzpicture}[use Hobby shortcut, scale=0.2, baseline=11pt]
	\foreach \k in {0,1.5,3} {
		\draw [white, double=black]([in angle=-90, out angle=90]1,\k) .. (0,1.5 + \k);
		\draw [white, double=black] ([in angle=-90, out angle=90]2,\k) .. (1,1.5 + \k);
		\draw [white, double=black, line width=1pt] ([in angle=-90, out angle=90]0,\k) .. (2,1.5 + \k); }
\end{tikzpicture}

%% file: tikz/kauff-dehn-3-proof-4.tex
\begin{tikzpicture}[use Hobby shortcut, scale=0.2, baseline=6pt]
	\foreach \k in {0,1.5} {
		\draw [white, double=black]([in angle=-90, out angle=90]1,\k) .. (0,1.5 + \k);
		\draw [white, double=black] ([in angle=-90, out angle=90]2,\k) .. (1,1.5 + \k);
		\draw [white, double=black, line width=1pt] ([in angle=-90, out angle=90]0,\k) .. (2,1.5 + \k); }
\end{tikzpicture}

%% file: tikz/kauff-dehn-3-proof-10.tex
\begin{tikzpicture}[use Hobby shortcut, scale=0.2, baseline=2pt]
	\draw [white, double=black]([in angle=-90, out angle=90]1,0) .. (0,1.5);
	\draw [white, double=black] ([in angle=-90, out angle=90]2,0) .. (1,1.5);
	\draw [white, double=black, line width=1pt] ([in angle=-90, out angle=90]0,0) .. (2,1.5);
\end{tikzpicture}

%% file: tikz/kauff-dehn-3-proof-12.tex
\begin{tikzpicture}[use Hobby shortcut, scale=0.2, baseline=3pt]
	\draw [white, double=black]([in angle=-90, out angle=90]1,0) .. (2,2);
	\draw [white, double=black, line width=1pt] ([in angle=90, out angle=90]0,0) .. (2,0);
	\draw [white, double=black] ([in angle=-90, out angle=-90]0,2) .. (1,2);
\end{tikzpicture}

%% file: tikz/kauff-dehn-3-proof-13.tex
\begin{tikzpicture}[use Hobby shortcut, scale=0.2, baseline=2pt]
	\draw [white, double=black]([in angle=-90, out angle=90]1,0) .. (0,1.5);
	\draw [white, double=black] ([in angle=-90, out angle=90]2,0) .. (2,1.5);
	\draw [white, double=black, line width=1pt] ([in angle=-90, out angle=90]0,0) .. (1,1.5);
\end{tikzpicture}

%% file: tikz/kauff-dehn-3-proof-14.tex
\begin{tikzpicture}[use Hobby shortcut, scale=0.2, baseline=3pt]
	\draw [white, double=black]([in angle=-90, out angle=90]1,0) .. (0,2);
	\draw [white, double=black, line width=1pt] ([in angle=90, out angle=90]0,0) .. (2,0);
	\draw [white, double=black] ([in angle=-90, out angle=-90]1,2) .. (2,2);
\end{tikzpicture}